\documentclass[10pt]{amsart}

\usepackage{amsmath}
\usepackage{enumerate}
\usepackage{hyperref}
\usepackage{amsopn}
\usepackage{amsmath, amsthm, amssymb, amscd, amsfonts, xcolor}
\usepackage[latin1]{inputenc}
\usepackage{graphicx}
\usepackage{amssymb}
\usepackage{pdfsync}
\usepackage{euscript}

\usepackage{xcolor}
\usepackage{enumitem}
\usepackage{multirow}

\input xy
\xyoption{all}

\usepackage[OT2,OT1]{fontenc}
\newcommand\cyr{%
\renewcommand\rmdefault{wncyr}%
\renewcommand\sfdefault{wncyss}%
\renewcommand\encodingdefault{OT2}%
\normalfont
\selectfont}
\DeclareTextFontCommand{\textcyr}{\cyr}

\theoremstyle{plain}

\newtheorem{teo}{Theorem}[section]
\newtheorem*{thm1}{Theorem A}
\newtheorem*{thm2}{Theorem B}

\newtheorem{cor}[teo]{Corollary}

\newtheorem{prop}[teo]{Proposition}
\newtheorem{lema}[teo]{Lemma}

\newtheorem*{prop1}{Proposition C}
\theoremstyle{definition}
\newtheorem{defi}[teo]{Definition}

\newtheorem{nota}[teo]{Remark}

\numberwithin{equation}{teo}

\def\gnul{E}
\def\lienul{\mathfrak{e}}
\makeatletter\renewcommand\theenumi{\@roman\c@enumi}\makeatother

\def\A{\mathbf{A}}
\def\ed{e}

\title{Homogeneous Riemannian manifolds with non-trivial nullity}

\author[A. J. Di Scala]{Antonio J. Di Scala}

\author[C.  Olmos]{Carlos   Olmos }
\thanks{A.J. Di Scala is member of GNSAGA of INdAM and of DISMA Dipartimento di Eccellenza MIUR 2018-2022, C. E. Olmos was supported by Famaf-UNC and CIEM-Conicet, and F. Vittone was supported by UNR and Conicet.}

\author[F. Vittone]{Francisco Vittone}


\begin {document}

\begin{abstract}
This is a substantially improved version of an earlier preprint of the authors
with a similar title.

We develop a general theory for irreducible homogeneous spaces $M= G/H$, in relation to
the nullity  distribution $\nu$ of their curvature tensor. We construct natural invariant
(different and increasing) distributions associated with the nullity, that give  a deep insight of such spaces. In particular, there must exist an order-two transvection,
not in the nullity, with null Jacobi operator. This fact was very important for finding
out the first homogeneous examples with non-trivial nullity, i.e. where the nullity distribution is not parallel.
Moreover, we construct irreducible examples of conullity $k=3$, the smallest possible, in any dimension. None of our examples admit a quotient of finite volume.
We also proved that $H$ is trivial and $G$ is solvable if $k=3$.

Another of our main results is that the leaves,  i.e. the integral manifolds, of the nullity are closed (we used a rather delicate argument). This implies that $M$ is a Euclidean affine bundle over the quotient by the leaves of $\nu$. Moreover, we prove that $\nu ^\perp$ defines a metric connection on this bundle with transitive holonomy or, equivalently,  $\nu ^\perp$ is completely non-integrable (this is not in general true for an arbitrary autoparallel and flat invariant distribution).

We also found some general obstruction for the existence of non-trivial nullity:
e.g., if $G$ is reductive (in particular, if $M$ is compact), or if
$G$ is two-step nilpotent.

\end{abstract}

\maketitle

\section {Introduction}

Given a Riemannian manifold $M$ with  curvature tensor $R$ and a point $p\in M$, the \textsl{nullity subspace} $\nu_p$ of $M$ at $p$ is defined as the subset of $T_pM$ consisting of those vectors that annihilate $R$, i.e., $$\nu_p=\{v\in T_p M\,:\, R_{\,\cdot,\,\cdot}\,v\equiv0\}.$$

The concept of nullity of the curvature tensor was first introduced by Chern and Kuiper in \cite{CK}. For a general Riemannian manifold, the dimension of the nullity subspace at a point $p$, called the \textsl{index of nullity} at $p$, might change from point to point. In the open and dense subset $\Omega$ of $M$ where the index of nullity is locally constant,  $q\mapsto \nu _q$  is an autoparallel distribution  with flat (totally geodesic) integral manifolds (the so-called leaves of the nullity).  Moreover, in the open subset of $\Omega$  where the index of nullity attains its minimum, the integral manifolds of $\nu$ are complete (cf. \cite{M}).
The distribution $\nu$ of $\Omega$  is called the \textit{nullity distribution} of $M$.  We  say that the nullity (distribution) $\nu$ is trivial  if it is the tangent space to a local de Rham flat factor of
$\Omega$, i.e.  $\nu$ is a parallel distribution with respect to the Levi-Civita connection or equivalently $M$ locally splits off the leaves (i.e. the integral manifolds) of the nullity distribution.   If $\nu = \{0\}$ we also say that $\nu$ is trivial  (in this case there must exist no local de Rham flat factor of $\Omega$).

\smallskip

Manifolds with positive index of nullity  have been studied by different authors  (see for example \cite{M}, \cite{M2}, \cite{G}, \cite{CM}, \cite{Ro1}, \cite{Ro2} and more recently, \cite {Sz},  \cite{BVK} and \cite{FZ}). The existence of non-trivial nullity has strong geometric implications. In particular, for the so-called CN2 manifolds, i.e., with codimension of the nullity at most two.

\medskip

 In this paper we study the nullity of  an irreducible homogeneous Riemannian manifold $M= G/H$.
  Since the curvature tensor $R$ is invariant under isometries, $M$ has constant index of nullity and so $\nu$ is a well defined  $G$-invariant autoparallel and flat  distribution of $M$.
 No progress had been made in this direction, except for some special cases where the nullity turned  out to be trivial, i.e. parallel (see, for instance, \cite {CFS}). Moreover, it was not known whether  an irreducible homogeneous space with non-trivial nullity could exist.
\smallskip

We develop a general  structure theory for such spaces, in relation to the nullity.
 Our approach is geometric and based on general facts about Killing fields, in contrast with the usual Lie algebra approach, i.e. we deal with Killing vector fields rather than with left-invariant objects. We use the so-called Kostant connection in the canonical bundle $E = TM \oplus \Lambda ^2 (TM)$, and his method  \cite{K} for computing the holonomy in terms of the Nomizu operators
 $\nabla X$, where $X \in \mathcal K (M)$ (the space of Killing fields).
\smallskip

A first question that naturally arises is whether the nullity foliation in a homogeneous Riemannian manifold is a homogeneous foliation, i.e. given by the orbits of an isometric group action. If $M$ has no Euclidean de Rham factor, we proved that this is never the case. Moreover, no Killing field $X\neq 0$ of
$M$ can be always tangent to the nullity, see Proposition \ref {KillTanNu}. Note that any single leaf of nullity is a homogeneous submanifold, but the subgroups of $I(M)$ that act on different leaves are in general only conjugate to each other.
\smallskip

 For dealing with the nullity $\nu$ we have to study non-trivial geometric distributions that are naturally associated to $\nu$ (when it is non-trivial). Namely, the osculating distributions
 $\nu ^{(1)}$, $\nu ^{(2)}$ of $\nu$, of first and second order.
 That is, if $\nu ^{(0)}:= \nu$, then $\nu ^{(i +1)}$ is obtained by adding to
 $\nu ^{(i)}$ the covariant derivative, in any direction, of fields that lie in $\nu ^{(i)}$ ($i= 0,1$). Since such distributions are $I(M)$-invariant one has, in particular, that
 $\nu ^{(1)} = \nu  + \hat \nu $, where $\hat \nu $ is the so-called adapted distribution of $\nu$. Namely, $\hat \nu _p$ is the  linear span of  $\{\nabla _{\nu _p}Z\}$, where $Z$ is a Killing field of $M$ induced by $G$.
 In our geometric construction it also appears another  natural $I(M)$-invariant distribution $\mathcal U$, the so-called bounded distribution, obtained by adding to $\nu$ the directions of the Killing fields whose normal component is bounded on a given leaf of $\nu$.  The distribution  $\mathcal U$ is $G$-invariant and $\nu$ is parallel in the directions of $\mathcal U$. Moreover, it is the largest $G$-invariant distribution with this property.
 \smallskip

 Our main results are summarized in the following theorems, which in particular relate the nullity to the so-called distribution of symmetry \cite {BOR}.  By making use of our construction, we were able to produce the first  examples of irreducible homogeneous manifolds with non-trivial nullity.  Moreover, we found, in any dimension, irreducible examples with codimension of the nullity equals to $3$, the smallest possible, and co-index of symmetry $2$, the smallest possible (see Theorem \ref{examples}).
 In, particular, in contrast to the compact case \cite {BOR}, there is no bound, in terms  of the co-index of symmetry, for the dimension of the space.
\smallskip

In any locally irreducible homogeneous Riemannian manifold, the codimension of the nullity must be always at least $3$. In fact, a CN2 space is semi-symmetric, and a semi-symmetric locally homogeneous space is locally symmetric due to a well known result of Z. I. Szab\' o (see \cite {Sz}, Proposition 5.1, and \cite {BVK}). Thus, a CN2 homogeneous space is locally the product of a flat factor and a surface of (non zero) constant curvature. We also obtain this result as a by-product of our main constructions.

%Thus, if the codimension is $2$ the nullity is parallel i.e. a CN2 homogeneous manifold splits off the leaves of the nullity.
\medskip

\

\begin {thm1}\label{thm1}
Let $M=G/H$ be a simply connected homogeneous Riemannian manifold  without a Euclidean de Rham factor. Assume that the nullity distribution $\nu$ is non-trivial and let  $k$ be its codimension. Let
 $\nu ^{(1)}, \nu ^{(2)}$, $\hat \nu$ and $\mathcal U$,  be the osculating,  of order $1$ and $2$, the adapted and the bounded  $G$-invariant distributions associated to $\nu$, respectively. Then

\vspace {.2cm}

 {\rm (1) } $\nu ^{(1)} = \nu + \hat \nu$ is autoparallel and flat, and $\mathcal U$ is integrable. Moreover,
we have the following inclusions (in particular, $k\geq 3$):
$$\{0\} \subsetneq \nu \subsetneq \nu ^{(1)}\subsetneq \nu ^{(2)}\subset \mathcal U
\subsetneq TM.$$
Moreover, the integral manifolds of $\nu$ and $\nu ^{(1)}$ are simply connected (and so isometric to a Euclidean space).
\vspace {.2cm}

 {\rm (2) }   For any $p\in M$, $v\in \hat \nu _p$,  there exists a transvection $Y\in \mathcal K ^G(M)$at $p$, with $Y_p=v$, and  the Jacobi operator $R_{\cdot , v}v$ is null. Moreover, $[Y,[Y, \mathcal K(M)]] = 0$ and $Y$ does not belong to the center of $\mathcal K ^G(M)$ (the Killing fields induced by $G$).  In particular, there exists such a transvection with $Y_p\notin \nu _p$.
\vspace {.2cm}

{\rm (3) }  The representation of the isotropy $H$ on $\nu _p^\perp$ is faithful
($p= [e]$). Moreover,
$\dim H \leq \frac 12 (k-2)(k-3)$ and so, if the equality holds and $G$ is connected,  $H\simeq \mathrm {SO}(k-2)$.  In particular,  if $k=3$, $H=\{e\}$.
\vspace {.2cm}

{\rm (4) }  If $k=3$, then  $G$ is solvable. Moreover, there exist irreducible examples in any dimension where $G$ is not unimodular, and so $M$ does not admit a finite volume quotient.

\vspace {.2cm}

\end {thm1}

\medskip

We should point out that we do not know any  example of an irreducible simply connected homogeneous Riemannian manifold $M=G/H$, with a non-trivial nullity,  such that  either $H$ is non-trivial or  $G$ is non-solvable ($G$ connected). There seems to be no geometric reason for the non-existence of such a  space.

\medskip

By making use of  some delicate arguments we were able to prove that the leaves of the  nullity foliation of an irreducible simply connected homogeneous space $M$  are closed (embedded) submanifolds.  Then we prove that $M$ is an affine bundle over the quotient of $M$ by the leaves of the nullity foliation with an affine connection. Moreover, $\nu ^\perp$ defines an affine connection on this bundle which has a transitive holonomy group, and so this distribution is completely non-integrable. The existence of a completely non-integrable geometric distribution on a Riemannian manifold has usually strong implications, as e.g. the so-called homogeneous slice theorem used in \cite {Th} and \cite {HL}.

\begin {thm2} Let $M=G/H$ be a simply connected homogeneous Riemannian manifold without a Euclidean de Rham factor, where $G=I(M)^o$.
Assume that $M$ has a non-trivial  nullity distribution $\nu$.
Then

{\rm (1) } Any integral manifold $N(q)$ of the nullity distribution $\nu$ is a closed  embedded submanifold of $M$ (or equivalently, the Lie subgroup $E^q$ of $G$ that leaves $N(q)$ invariant is closed).

{\rm (2) } $M$ is  the total space of a Euclidean  affine bundle over the quotient $B=G/E^p$ of $M$ by the leaves of nullity  foliation  with standard fiber $N(p)\simeq \mathbb R^{\mu}$ ($p=[e]$).

{\rm (3) } $\nu ^\perp$ defines an affine metric connection, in the sense of Ehresmann, on the affine bundle
$M\to B$. Moreover, the holonomy group  associated to $\nu ^\perp$ is transitive  (or equivalently,
$\nu ^\perp$ is completely non-integrable).

\end {thm2}

\smallskip

Observe that  part (3) above is not true for an arbitrary autoparallel and flat $G$-invariant distribution of a homogeneous Riemannian manifold $M=G/H$.
In fact, in a symmetric space of the non-compact type, presented as a solvable Lie group $M=G$, the normal spaces of a  foliation by parallel (geometric) horospheres define an autoparallel and flat $G$-invariant distribution $\mathcal A$, but $\mathcal A^\perp$ is
integrable.

\medskip

By making use of our main results, we  found some  obstructions for the existence of non-trivial nullity.  In fact, the  existence of transvections of order $2$ imposes general restrictions on the presentation group  $G$. Namely,

\begin {prop1}\label {ProB} Let $M=G/H$ be a simply connected homogeneous Riemannian manifold  without a Euclidean de Rham factor. Then the nullity distribution is  trivial, with any $G$-invariant metric, in any of the following cases:
\vspace {.2cm}

\rm (a) If the Lie algebra of $G$ is reductive,  i.e. the direct sum of an abelian ideal and a semisimple ideal  (in particular, if $M$ is compact).
\vspace {.2cm}

{\rm (b)} If  the Lie algebra $\mathfrak g$ of $G$ is $2$-step nilpotent.

\vspace {.2cm}

\end {prop1}

\smallskip

\noindent
{\bf Open questions.} Let $M=G/H$ be a simply connected irreducible Riemannian
homogeneous manifold with positive index of nullity.
\begin {itemize}
 \item [{\bf - }] Is the isotropy group $H$ always trivial?

 \item [{\bf - }] Is $G$ necessarily solvable?
 \end {itemize}

\smallskip

The reader interested in the proofs of Theorems A, B and Proposition C  and that would like to avoid  the preliminaries can find them in Section \ref{demostraciones}.
\vspace {.2cm}

{\small This is a substantially improved new version of an older preprint of the authors
that circulated with a similar title.}\\
\section {Preliminaries and basic facts}\label {Preliminaries}

Let $(M,\langle \, ,\, \rangle)$ be a (connected) complete Riemannian manifold with Levi-Civita connection $\nabla$. A vector field $X$ of $M$ is called a Killing field if
\begin{equation}\label{killingcond}
v\mapsto \nabla _v X
 \end{equation} is a skew-symmetric endomorphisms  of $T_pM$, for all $p\in M$.
 Such a condition is called the {\it Killing equation} and  reflects  the fact that the  flow of $X$ is by isometries.

Let
$ I(M)$ denote the  Lie group of isometries of $M$. The Lie algebra
$\mathrm {Lie}(I(M))$ of $ I(M)$ is naturally identified with
the Lie algebra $\mathcal K (M)$ of Killing fields of $M$.
Namely, the map $z\overset {j}{\mapsto} z^*$ is a linear isomorphism
from $\mathrm {Lie}(I(M))$ onto $\mathcal K (M)$ that satisfies
\begin{equation}\label{laspe2}[x,y]^*= - [x^*, y^*],
\end{equation}
 where
$$z^*_q = z.q := \frac {\text {\rm d}}{\text {\rm d}t}_{
\vert 0} \text {\rm Exp}(tz)q.$$
 In fact, let
$f: I(M) \to M$ be the map $f(g) = g(p)$, $p\in M$ fixed. Then the right
invariant vector field with initial condition $z\in T_e I(M)$ is $f$-related to the Killing field $z^*$. The vector field $ z^*$ is called
the Killing field {\it induced} by $z\in \mathrm {Lie}(I(M))$.

\begin{nota}\label{LeftRight}
If one should define the Lie algebra  $\mathrm {Lie}(I(M))$  by using right-invariant
vector fields instead of left-invariant vector fields, then the map $j$
would be a Lie algebra isomorphism
(see, for instance,  A.2 in  \cite {BCO}).
\end{nota}

If $G$ acts by isometries on $M$ and $z\in\mathfrak{g}=\mathrm{Lie}(G)$, then the field $z^*$ is called a  Killing field of $M$ induced by $G$. We will denote the set of such vector fields by $\mathcal K ^G(M)$. If the action of  $G$  on $M$ is not effective, there could exist non-zero elements $z\in\mathfrak{g}$,  such that the corresponding $z^*\equiv 0$.

\medskip

Let $X\in \mathcal K (M)$. The initial conditions of $X$ at  $p\in M$ are given by the pair
$$(X)^p: =(X_p, (\nabla X)_p)\in T_pM\oplus \Lambda^2(T_pM).$$
where $(\nabla X)_p$ denotes the skew-symmetric endomorphism defined by equation (\ref{killingcond}).
  These conditions completely determine the Killing field $X$, in the sense that two Killing fields with the same initial conditions at any fixed point $p$  must coincide on $M$.

A Killing field $X$, besides the Killing equation, satisfies
the following identity, for all $p\in M$, $u,v\in T_pM$
\begin {equation}
\label {affK}
\ \ \ \ \ \ \ \ \ \ \ \ \nabla ^2_{u,v}X = R_{u,X_p}v
\ \ \ \ \ \ \ \ \ \ \mathit {affine \ Killing \ equation}.
\end {equation}
The affine Killing equation reflects the fact the flow of $X$ preserves the Levi-Civita connection.

Equations (\ref{killingcond}) and (\ref{affK}) motivate the introduction of the so-called \textsl{Kostant connection}  $\tilde{\nabla}$ on the vector bundle $$E := TM \oplus \Lambda  ^2(TM)$$ (see \cite {K, CO}). Here $\Lambda  ^2(T_pM)$ is, as usual, identified with the skew-symmetric  endomorphisms of $T_pM$.
The bundle $E$ is called the {\it canonical bundle} and  $\tilde \nabla$ is given by
\begin {equation}\label {Konstant}
\tilde \nabla _u (Z, B) = (\nabla _u Z - Bu, \nabla _u B - R_{u, Z_p})
\end {equation}
 $u\in T_pM$, where $(Z, B)$ is a section of $E$ and $R$ is the curvature tensor of $M$.
The Killing  fields of $M$ are naturally identified with the parallel sections
of $E$ in the following way: $(X, B)$ is a parallel section of $E$ if and only if $X$ is a Killing field of $M$ and $B= \nabla X$.

If $X$ is a Killing field, then the section $q\mapsto (X_q,(\nabla X)_q)$ is called the {\it canonical lift} of $X$ to $E$.

The Kostant connection allows us to determine the initial conditions of a Killing field $X$ at any $q\in M$ if we know the initial conditions
$(X) ^p$ at a fixed $p$. In fact, we must compute the parallel transport, in the Kostant connection, of $(X) ^p$ along any curve from $p$ to $q$ (in particular, by using a geodesic).

From the affine Killing equation and the Bianchi identity one can determine the initial conditions at $p$ of the bracket $[X,X']$ of any two Killing fields
in terms of the initial conditions $(X)^p = (v, B)$, $(X')^p = (v', B')$
 (see  Lemma 2.4 of \cite {R}). Namely,
\begin{equation} \label{initialbracket}
([X,X'])^p = (B'v-Bv', R_{v, v'} - [B,B']).
\end{equation}

This equation gives a useful formula for computing the curvature in terms of Killing fields $X$ and $Y$:
\begin {equation} \label {curvatureformula}
 R_{X_p, Y_p} = (\nabla [X,Y])_p + [(\nabla X)_p , (\nabla Y)_p].
\end {equation}

The well-known Koszul formula gives the Levi-Civita connection $\nabla$ in terms of
brackets of vector fields and scalar products.
Since the Lie derivative of the metric tensor along any Killing vector field is zero, we have the following expression for $\nabla$ in terms of Killing fields $X, Y, Z$ (see (3.4) in p. 617 of \cite {ORT})
\begin {equation} \label {fundamentalequation}
2\langle \nabla _XY , Z\rangle = \langle [X,Y],Z \rangle
+ \langle [X,Z],Y \rangle + \langle [Y, Z], X \rangle.
 \end {equation}

 \begin {nota} \label {limitK} Regarding Killing fields as sections of the canonical bundle $E$, one can easily  prove the following fact: let $Z^n$ be a   sequence of Killing fields  on $M$ induced by $G$ such that, for some $p\in M$, their initial conditions at $p$, $(Z^n)^p = (Z^n_p, (\nabla Z^n)_p)$ converge to $(v,B)$. Then $(v,B)$ is the initial condition at $p$ of a Killing field $Y$ induced by $G$. Moreover, $((Z^n)_q, (\nabla Z^n)_q)\to (Y_q, (\nabla Y)_q)$, for all $q\in M$.
\end {nota}

\smallskip

 \subsection {Parallel transport along integral curves of Killing fields}
 \label {parallelintegral} \
 \smallskip

 Let $X$ be a Killing field and let $\phi _t$ be its associated flow. Observe that such a flow is always of the form $q\mapsto \mathrm
 {Exp}(tz)q$, for some $z$ in the Lie algebra of the isometry group.
 Let $p\in M$ and let $c(t)= \phi _t (p)$ be the integral curve of $X$ by $p$. Let $\tau _t$ denote the parallel transport along $c(t)$, form $0$ to $t$. Then it is  not difficult to show that
 $\tau _{t}^{-1}\circ \mathrm {d}\phi _t : T_pM \to T_pM$ is a $1$-parameter subgroup of
 linear isometries. Moreover, (see e.g. Remark 2.3 of \cite {OS}),
 \begin {equation} \label {etnabla}
 \tau _{t}^{-1}\circ \mathrm {d}_p\phi _t  =\mathrm{e} ^ {t(\nabla X)_p}.
 \end {equation}
 In fact, if $\gamma ' (0) =v\in T_pM$,
\begin{eqnarray*}
\tfrac {\, \mathrm d}{\mathrm d t}_{\vert 0}
  \tau _{t}^{-1}\circ \mathrm {d}_p\phi _{t}  (v) &=&
 \tfrac {\, \mathrm {D} }{\partial t}_{\vert 0}
 \tfrac {\, \partial }{\partial s}_{\vert 0}
 \phi _t (\gamma (s))  = \tfrac {\, \mathrm {D} }{\partial s}_{\vert 0}
 \tfrac {\, \partial }{\partial t}_{\vert 0}
 \phi _t (\gamma (s)) \\
 & =&
 \tfrac {\, \mathrm {D} }{\mathrm d s}_{\vert 0}
 X_{\gamma (s)}= \nabla _vX.
\end{eqnarray*}

 \begin {nota}\label {accumulation}
 If $H$ is a $1$-dimensional Lie subgroup of a compact Lie group
 $K$, then the closure $T$ of $H$ is an abelian compact subgroup of $K$, i.e. a torus. From this it is not hard to see that there is a sequence of real numbers  $\{t_n\}_{n\in \mathbb N}$, which tends to $+ \infty$ and such that
 $$\tau _{t_n}^{-1}\circ \mathrm {d}_p\phi _{t_n}  =\mathrm{e} ^ {t_n(\nabla X)_p}$$
 tends to the identity transformation of $T_pM$ (or to any other element
  of the closure of $\{\mathrm{e} ^ {t(\nabla X)_p}: t \in \mathbb R\}$).
\end {nota}

\begin {nota}\label {isotropy} Let $X$ be a Killing field  that belongs to the isotropy algebra at $p$, i.e., $X_p =0$. Let, in the previous notation,  $\phi _t$ be the flow associated to $X$ and $c(t) = \phi _t (p) \equiv p$. Then from
  (\ref{etnabla}) it follows that $\mathrm {d}_p\phi _t = \mathrm {e}^{t(\nabla X)_p}$. So, via the isotropy representation at $p$, $\phi _t $ is identified with  $\mathrm {e}^{t(\nabla X)_p}$.
  \end {nota}

\smallskip

\subsection {Holonomy of homogeneous spaces:  Kostant's results} \label {Kostantmethod}\
\smallskip

Let $M= G/H$ be a homogeneous Riemannian manifold and let $\mathcal K ^G(M)$ be the space of Killing fields on $M$ induced by $G$. Let $p\in M$ and let
 $\tilde {\mathfrak h} (p)$ be the Lie subalgebra of $\mathfrak {so}(T_pM)$ which is algebraically spanned by the set $\{(\nabla X)_p:X\in \mathcal K ^G(M)\}$. Kostant proved in \cite {K} that $\tilde {\mathfrak h} (p)$ contains the holonomy algebra $\mathfrak {hol}(p)$ of $M$ at $p$  and it is contained in its normalizer $\mathfrak n (p)$ of $\mathfrak {hol}(p)$ in $\mathfrak {so}(T_pM)$, i.e,
\begin{equation}
\mathfrak {hol}(p)\subset \tilde {\mathfrak h} (p)\subset \mathfrak n (p).
\label{holalg}
\end{equation}
 Moreover, if $M$ is locally irreducible    and it is not   Ricci flat, he proved that
$\tilde {\mathfrak h} (p) $ coincides with $\mathfrak {hol}(p)$ (for a modern treatment of this subject see  the survey \cite {CDO}). Since Alekseevskii and Kimelfeld proved in \cite{AK} that a homogeneous and Ricci flat space must be flat, one has that
$\tilde {\mathfrak h} (p) = \mathfrak {hol} (p)$
for a locally irreducible homogeneous Riemannian manifold.

On the other hand, for any locally irreducible Riemannian  manifold, the normalizer of the holonomy algebra, inside the orthogonal algebra, properly contains the holonomy algebra if and only if the space is K\"ahler and Ricci flat (see e.g. \cite {CDO}, \cite [Prop. 5.2.3] {BCO}).
This, together with Kostant result (\ref{holalg}), implies that for a possibly reducible  homogeneous space $M$,
\begin{equation}\label{deRham1} \mathfrak {hol}(p) = \tilde {\mathfrak h} (p) \end{equation}
if $M$ has no (local) Euclidean de Rham factor or $M$ has a Euclidean de Rham factor of dimension $1$.

\begin {cor}\label {KostantdeRham}
Let   $M=G/H$ be a homogeneous Riemannian manifold. Assume that   there is a non-trivial subspace $\mathbb V$ of $T_pM$ which is invariant by $(\nabla X)_p$, for all
$X\in \mathcal K ^G(M)$. Then   $\mathbb V$ extends locally to a parallel distribution of $M$ and so $M$ locally splits.
\end {cor}

\begin {proof}
The subspace $\mathbb V$ is invariant by $\tilde{\mathfrak {h}}(p)$,
and therefore by $\mathfrak {hol} (p)$ from (\ref{holalg}). Then, by Remark \ref{isotropy}, $\mathbb V$ locally extends to a parallel non-trivial distribution. So, de  Rham decomposition theorem applies.
But, for the sake of self-completeness,  let us do a direct proof.

 From Remark \ref {isotropy} $\mathbb V$ is invariant under the isotropy algebra. Since we are working locally we may assume that $H$ is connected. Then $\mathbb V$ extends to a $G$-invariant distribution $\mathcal D$ on $M$. Moreover, for any $q\in M$, $\nabla _{\mathcal D_q}X \subset \mathcal D_q$, for all $X\in \mathcal K ^G(M)$.
Let $\xi$ be a field on $M$ that lies in $\mathcal D$ and let
$X\in \mathcal K ^G(M)$ be arbitrary. Since $\mathcal D$ is $G$-invariant then
$[X,\xi]$ lies in $\mathcal D$. But $[X,\xi] = \nabla _X\xi  -\nabla _{\xi} X$. Since
$\nabla _{\xi} X$ lies in $\mathcal D$, then $\nabla _X\xi $ lies in $\mathcal D$.
Then $\mathcal D$ is a non-trivial parallel distribution and $M$ splits locally.
\end {proof}

\smallskip

\subsection {The index of symmetry}\label {index}\
\smallskip

Let $M$ be a homogeneous Riemannian manifold. A Killing field $X$ on $M$ is called a
{\it transvection} at $q$ if $$(\nabla X)_q=0.$$
If $X$ is a transvection at $q$ it follows, from (\ref {etnabla}),  that
 $\gamma (t) :=  \phi _t (q)$ is a geodesic in $M$ and $\mathrm d_q \phi _t$ gives the parallel transport along $\gamma (t)$.

We now introduce some basic definitions that were given in \cite {ORT} (see also \cite {BOR}).

The {\it  Cartan subspace} at $q$ is
\begin {equation}\label {Cs}
\mathfrak p ^q := \{X\in \mathcal K (M): X \text { is a transvection at }q\},
\end {equation}
 the {\it symmetric isotropy algebra} at $q$ is
$$\mathfrak t  ^q := [\mathfrak p ^q, \mathfrak p ^q],$$
 and the {\it symmetric subspace} at $q$ is
$$\mathfrak s_q : = \mathfrak p ^q.q =
\{X_q\ :\, X\in \mathcal K (M),\, X \text { is a transvection at }q\}. $$

It turns out that
$$\tilde {\mathfrak  g}^q = \mathfrak t  ^q \oplus \mathfrak p ^q$$
is an involutive Lie algebra, the so-called {\it Cartan algebra} at $q$.
(We have used the notation $\tilde {\mathfrak g}^q$ instead of the more natural
$ \mathfrak{g}^q$, as in the references, in order to be consistent  with further notation).

Since we are assuming that $M$ is homogeneous,  all the previous objects are conjugate to each other by an isometry if we change the base point. In this way $\mathfrak s$ defines an $I(M)$-invariant  distribution on $M$ which is autoparallel. So it is well defined the so-called
{\it index of symmetry}, $i _{\mathfrak s}(M)$,   as the dimension over $M$ of the distribution
$\mathfrak s$.

The integral manifold $L(q)$ of $\mathfrak s$ through $q\in M$ is a totally geodesic submanifold of $M$, called the {\it leaf of symmetry} through $q$.
The leaves of symmetry are globally symmetric spaces as it follows from Corollary 2.3 in \cite {BOR}.

The autoparallel subdistribution of $\mathfrak s$ associated to the flat local de Rham factors of the leaves of symmetry will be denoted by $\mathfrak s ^0$. The set of associated transvections at $q$ will be denoted by $\mathfrak p^q_0$, i.e, $\mathfrak p^q_0=\{X\in \mathfrak{p}^q\,:\,X_q\in \mathfrak{s}^0_q\}$.
Observe that $\mathfrak p^q_0$ is the abelian part of $\mathfrak p^q$. That is to say,
\begin{equation}\label{flatfactor} \mathfrak p^q_0 = \{ X \in {\mathfrak p}^q \mbox{ : } [X,{\mathfrak p}^q] = 0  \} \, .
\end{equation}

The {\it co-index of symmetry} is the codimension of the distribution of symmetry.
If $M$  is  not locally symmetric,   then its co-index of symmetry
is at least $2$. This was shown in \cite{BOR} for the compact case and by Reggiani
 for the general case in \cite[Theorem 2.2]{R}.

\begin{nota}
Observe that all the previous geometric objects have been defined using all the Killing fields of $M$, i.e, Killing fields induced by the whole isometry group $I(M)$. In general, if $M=G/H$,  a transvection at $q$ may not be  a Killing field induced by the presentation group $G$   if $G$ is a proper  subgroup of
$I(M)^o$. For instance,  if  $S^3 = \mathrm {Spin_3}$, then no transvection of the sphere is a Killing field induced by the presentation group. 
\end{nota}

We now generalize Lemma 3.3 in \cite{ORT} for the case where $M$ is not necessarily compact.
\begin{lema}\label {EFF}
Keeping the notations of this section, the Lie subgroup $\tilde{G}^q\subset I(M)$ whose Lie algebra is $\tilde {\mathfrak g} ^q$, acts almost effectively on $L(q)$.
\end{lema}

 \begin{proof}
 We are going to use the following fact: for a Riemannian homogeneous space $M$, the Killing form $B$ of $\mathrm{Lie}(I(M))$  is negative definite when restricted to  $\mathrm{Lie}(I(M)_q)$
 (the Lie algebra of the isotropy at $q$).

Consider the ideal $\mathfrak h$ of $\tilde{\mathfrak{g}}^q$ given by the elements $X$ such that $X_{|L(q)}\equiv 0$. Then $\mathfrak{h}\subset \mathfrak{t}^q$. Since $[\mathfrak{t}^q,\mathfrak{p}^q]\subset \mathfrak{p}^q$ we get that $$[\mathfrak{h},\mathfrak{p}^q]=0.$$
Therefore, if $X\in \mathfrak{h}$ and $Y,Z\in \mathfrak{p}^q$, then
$$B(X,[Y,Z])=B([X,Y],Z)=0.$$
So $B(X,\mathfrak{t}^q)=0$, and then $X\equiv 0$.
 \end{proof}

\medskip

\subsection {The nullity of the curvature tensor} \label {nullity}\
\smallskip

Let $M$ be a Riemannian manifold. The nullity of the curvature tensor $R$ at $p\in M$ is
$$\nu _p := \{v \in T_pM: R_{v, x} = 0, \forall x \in T_pM \}$$
or, equivalently, due to the identities of the curvature tensor,
$$\ \ \ \ \ \nu _p := \{v \in T_pM: R_{x, y}v = 0, \forall x, y \in T_pM\}.$$

The nullity $\nu$ defines a (differentiable) distribution in the open and dense subset $\Omega$ of
$M$ where the dimension $\dim (\nu _q)$ is locally constant. Moreover, as it is well-known, it is an autoparallel distribution (or equivalently, it is integrable with totally geodesic integral manifolds). For the sake of self-completeness let us show this fact. Let $X,Y,Z, W$ be arbitrary vector fields in a connected component of $\Omega$ such that $X$ and $Y$ lie in $\nu$.
Then, a direct calculation shows that
$(\nabla _ZR)_{W,X}Y = 0 = (\nabla _WR)_{X,Z}Y$. Then, by the second Bianchi identity,  one has that
\begin {eqnarray*}
0 &=& (\nabla _XR)_{Z, W}Y\\
&=& \nabla _XR_{Z, W}Y - R_{\nabla _XZ, W}Y - R_{Z, \nabla _XW}Y -
R_{Z,W}\nabla _XY \\
&=& -R_{Z,W}\nabla _XY
\end {eqnarray*}
which shows that $\nabla _XY$ lies in $\nu$.

\begin {lema}\label {ABC} Let $M$ be a Riemannian manifold.
Let  $\gamma _v (t) $ be  a geodesic everywhere tangent to $\nu$,  with $\gamma (0) =p$, $\gamma' (0) = v$. Denote by $\tau_t$ the parallel transport along $\gamma(t)$ from $0$ to $t$.  Let $X$ be an arbitrary Killing field on $M$.
Then

\begin {itemize}
\item [(i)] $X_{\gamma (t)} =  \tau _t (X_p) + t \tau _t(\nabla _v X)$.

\item [(ii)] $\nabla _{\gamma'(t)} (\nabla X) =0$, i.e.,  $\nabla X$ is parallel along $\gamma (t)$, or equivalently
$$(\nabla X)_{\gamma (t)}= \tau _t ((\nabla X)_p) : = \tau _t
\circ (\nabla X)_p \circ \tau _t^{-1}.$$
\end {itemize}
\end {lema}

\begin {proof} Since $X$ is a Killing field then $X_{\gamma (t)}$ is a Jacobi field along
$\gamma (t)$. Observe that the Jacobi operator $R_{\cdot , \gamma' (t)}\gamma'(t) =0$.
Then the Jacobi equation yields  $\frac {\mathrm {D}^2}{\mathrm d t} X_{\gamma (t)}= 0$. This shows (i).

 Part
(ii) follows immediately from formula (\ref {affK}).
\end {proof}

\smallskip

Let now  $M$ be a homogeneous Riemannian manifold. Then the nullity distribution $\nu$,
being a geometric object,  is invariant under any isometry $g\in I(M)$, i.e., $g_*(\nu)
= \nu$. Thus,  $\dim (\nu_q)$ does not depend on $q\in M$ and therefore $\nu $ is a (smooth) distribution in $M$.

We will denote by $N(x)$ the so-called {\it leaf of nullity} by $x$, i.e, the totally geodesic (maximal) integral manifold of $\nu$ that contains
$x$. Then if $g\in I(M)$, we must have that  $gN(x) = N(gx)$,
for all $x\in M$.

\smallskip

 As a first consequence of the results of Section \ref{Kostantmethod} one can give a simple proof of Proposition C for the case where $M$ is compact:

\begin {prop}\label {compact2}

Let $M$ be a compact  homogeneous Riemannian manifold, which does not split off, locally, a flat de Rham factor. Then the distribution of nullity is trivial.
\end {prop}

\begin {proof}

Let $\gamma _v (t)$ be a non-constant geodesic tangent to the nullity distribution and let $X \in \mathcal K (M)$ be arbitrary. Since $M$ is compact, then $X$ is bounded and so, by  Lemma \ref {ABC},
$\nabla _vX=0$. By Corollary \ref {KostantdeRham}, since $M$ is homogeneous, $M$ splits locally the direction of $v$. A contradiction.
\end {proof}

\subsection{Homogeneous flat spaces} \label {homogeneous flat spaces}\
\smallskip

The main purpose of this section is to prove the well-known result that any connected transitive Lie subgroup $G$ of the isometry group of $\mathbb{R}^n$ must contain a pure translation. Moreover we will prove that any isometry of $\mathbb{R}^n$ that commutes with all the elements of such a group $G$ must be a translation.  As it will become clear in sections \ref{homogeo}, \ref{sec61} and \ref{LTA}  (the following sections), the constructions and techniques presented in the proofs of the following lemmas (here) can be easily adapted to homogeneous spaces with non-trivial nullity.

\smallskip

%\begin {nota}  \label {homogeneous flat spaces}
\begin{lema} \label{transtrasl}
Let $G$ be a  connected Lie subgroup of  the full isometry group $I(\mathbb R^n)= \mathrm {O}(n)\ltimes \mathbb R^n$. Assume that $G$ acts transitively on $\mathbb R^n$. Then $G$ must contain a pure translation.
\end{lema}
 \begin{proof}
 We identify any element $X$ of the Lie algebra $\mathfrak  {so}(n)\ltimes \mathbb R^n$ as a Killing field of the Euclidean space. Namely, if $X = B +v$, $B\in \mathfrak {so}(n)$ and $v\in \mathbb R^n$, then
$X_q = Bq +v$. Observe that for any $q\in \mathbb R^n$ and any $w\in T_q\mathbb R^n$, $$\nabla_wX = Bw,$$ and so $\nabla X$ is a parallel  skew-symmetric tensor. Therefore, $X$ is a transvection if and only if $B=0$, i.e., $X$ is a pure translation.

We shall see that there always exists a non-trivial transvection in the Lie algebra $\mathfrak g$ of $G$.
In fact, let $X = B +v \in \mathfrak g$,  $B\neq 0$ and let $w\in \mathbb R^n$ such that $Bw\neq 0$.
For each $t \in \mathbb{R}$ put $q_t = tw$. Then $\Vert X_{ q_t}\Vert \to \infty$, if $t\to +\infty$. Since $G$ acts transitively on $\mathbb R^n$, for each $t\in \mathbb{R}$ we can choose an element $g_t \in G$ such that  $g_t(q_t) = 0$. Set $X^t = (g_t)_* (X)\in \mathfrak{g}$.
Then $X^t_0 =  t\mathrm {d}g_t (Bw) +\mathrm {d}g_t (v)$ and $\nabla X^t =
\mathrm {d}g_tB \mathrm {d}(g_t)^{-1}$. Let now $t_k\to +\infty$ be such that the sequence (of constant norm)
$\mathrm {d}g_{t_k} (Bw)$ converges to some $0\neq u\in \mathbb R^n$. Then
$\frac { 1\, }{t_k}X^{t_k}$ converges to the transvection associated to $u$.
\end{proof}

 We will now introduce some technical results and fix some notation that we will need for the proof of our last lemma.  Let $G$ be a transitive Lie group of isometries of $\mathbb{R}^n$ with Lie algebra $\mathfrak{g}$ and let $\mathfrak a\neq \{0\}$ be the ideal of $\mathfrak g$ that consists of all the transvections in $\mathfrak g$. Let $$\mathbb V = \{X_0:X\in \mathfrak a\}.$$ One has, from the fact that
$\mathfrak a$ is an ideal of $\mathfrak g$, that
$B \mathbb V \subset \mathbb V$, where  $X= B +v \in \mathfrak g$ is arbitrary. Therefore, one also has that $B\mathbb{V}^{\perp}\subset \mathbb{V}^{\perp}$.

Consider the $G$-invariant distribution $\mathcal{D} ^\perp$ on $\mathbb R ^ n$ defined by $\mathbb V^\perp$ such that $\mathcal D^\perp _0 = \mathbb V^\perp $. Let us identify $\mathbb V^\perp$ with the integral manifold of $\mathcal D^\perp$ by $0$.
Since $G$ acts transitively on $\mathbb R^n$,  for any $w\in \mathbb V^\perp$ there exists $X \in \mathfrak g$ of the form $X = B + w$, where $B\neq 0$. Then $X_0\in \mathbb{V}^{\perp}$ and so $X_{|\mathbb{V}^\perp}$ is always tangent to $\mathbb{V}^\perp$. Therefore, the Lie subgroup $G'$ of $G$ that leaves $\mathbb V^\perp$ invariant is transitive on this subspace.

Let $X=B+v\in \mathfrak{g}$ be arbitrary. Let us see that the restriction $B_{\vert \mathbb V^\perp} =0$. In fact, if $B_{\vert \mathbb V^\perp} \neq 0$, the same limit argument used in the proof of Lemma \ref{transtrasl} for
$g_t\in G'$ would lead to the construction of a non-trivial transvection of $\mathbb R^n$ in a direction perpendicular to $\mathbb{V}$. This is a contradiction that proves that $$B_{\vert \mathbb V^\perp} =0.$$

Observe that all the integral curves of $X\in \mathfrak a$ are geodesics (i.e.,  lines) that lie in the $G$-invariant distribution $\mathcal D$ defined by $\mathbb V$. If $X\in \mathfrak g$ is such that $X_0\in  \mathbb V^\perp$, then the integral line of $X$ with initial condition $0$ is also a geodesic.

\begin{lema} \label{isoconmuta}
Let $G$ be a  connected Lie subgroup of  $I(\mathbb R^n)$ which  acts transitively on $\mathbb R^n$. Let $g\in I(\mathbb R^n)$ be an isometry that commutes with all the elements of $G$. Then $g$ is a translation.
\end{lema}

\begin{proof}
 We will keep the notations of the previous paragraphs. Observe first that $g\in I(\mathbb R^n)$ commutes with all the elements of $G$ if and only if $g_* (X) = X$ for all $X =B + z \in \mathfrak g$, the Lie algebra of $G$.

Let us write $g(x) = Ux + u$, where $U\in \mathrm O(n)$, $u\in \mathbb R^n$. Then, taking $X \in \mathfrak a$, one has that $U$ leaves $\mathbb V$ invariant and $U_{\vert \mathbb V}$ is the identity. Then, for an arbitrary $X = B + z$, one has that $UBU^{-1} = B$, since
$B_{\vert \mathbb V^\perp} = 0$. Let us assume that $z\in \mathbb V ^\perp$.  Then
\begin {eqnarray} \label {EU} g_* (X)_q &=& \mathrm {d}g (X_{g^{-1}q}) =
U (B ( U^{-1}(q-u))+ z) \\
& =& B(q-u) + Uz  = X_q = Bq +z\nonumber
\end {eqnarray}
and so $$(U-Id)z = Bu.$$
But the left hand side belongs to $\mathbb V ^\perp$ and the right one to $\mathbb V$. Then $Uz = z$, for all $z\in \mathbb V^\perp$, since there is a Killing field induced by $G$ in any direction.  Then $U= Id$, since $U$ acts trivially on $\mathbb V$, and $g$ is a translation.
\end{proof}

 As a consequence of the results of this section one obtains the well-known

\noindent
\textbf{Fact:}
any homogeneous  flat Riemannian manifold is a product of a torus and a Euclidean space (cf. \cite{AK}).

Let $\Gamma \subset I(\mathbb R^n)$ be a discrete subgroup that acts properly discontinuously and commutes with $G$.  Then  from Lemma \ref{isoconmuta}  $\Gamma$ consists of translations. So the quotient space of $\mathbb R^n$ by $\Gamma$ is a $G$-homogeneous space which is the Riemannian  product of a flat torus by a Euclidean space.
This implies that any flat homogeneous Riemannian manifold is such a product and so it is isometric to an abelian Lie group with an invariant metric  (the presentation group $G$ may be non-abelian).
% So, any homogeneous  flat Riemannian manifold is a product of a torus and a Euclidean space.

\medskip

Before concluding this section, let us recall the well-known fact that a Lie group of isometries that acts simply transitively on $\mathbb R^n$ must be ($2$-step) solvable (see e.g. \cite{Al}). In fact, this can also be obtained using the techniques developed in this section.
%\end {nota}

\section {The nullity of  homogeneous spaces} \label{secnul}

\subsection {The osculating  distributions of the nullity}\label {der}\
\smallskip

Let $M=G/H$ be a  Riemannian homogeneous manifold and assume that its nullity distribution $\nu$ is non-trivial. We also assume that $\nu$ is not parallel. Otherwise, $M$ would split off, locally, a flat factor.

Let us consider the {\it osculating distribution} $\nu^{(1)}$ associated to the nullity distribution. Namely, if  $C^\infty (\nu)$ are the tangent fields of $M$ that lie in $\nu$,

$$\nu ^{ (1)}_q = \nu _q +  \mathrm{span}\,  \{ \nabla _w X : X \in C^\infty (\nu), w\in T_qM\}.$$

Then $\nu ^{ (1)}$ is a $G$-invariant distribution that properly contains $\nu$, since
$\nu$ is non-parallel.  It is not hard to see that one only needs to consider $\nabla _w X$, for $X$ in some family of  fields  that lie in $\nu$  and such that $X_q$ span $\nu _q$.

The osculating  distribution can be defined for any $G$-invariant distribution $\mathcal H$ and $\mathcal H ^{(1)}$ is also a $G$-invariant distribution. The
osculating distribution of order $k$ is defined as
$$\mathcal H ^{k} = (\mathcal H ^{k-1})^{(1)}, $$
where $\mathcal H ^{(0)}:= \mathcal H$.

\begin {lema}\label {derived} Let $\mathcal H$ be a  $G$-invariant distribution
of $M$. Then
  $$ \mathcal H ^{ (1)}_q = \mathcal H _q +  \mathrm {span}\,  \{\nabla _v Z: Z\in \mathcal K ^G (M) , v\in \mathcal H_q\}.$$
\end {lema}

\begin {proof} Since  $\mathcal H$ is $G$-invariant, the flow of  any Killing field $Z$, induced by $G$,  leaves $\mathcal H$ invariant. Then, differentiating   the flow, one obtains  that
$[Z,X]$ lies in $\mathcal H$, if $X$ lies in $\mathcal H$. So, if $X$ lies in $\mathcal H$,
$$\nabla _{Z_q}X = \nabla _{X_q}Z + [Z, X]_q \sim \nabla _{X_q}Z  \ (\text {mod }
\mathcal H _q).$$
Since $M$ is homogeneous, there are Killing fields in any arbitrary direction. This proves the lemma.
\end {proof}

 We will show later that
 $\nu ^{(1)}$ is an  autoparallel and flat distribution (that properly contains $\nu$). Moreover, we will show  that $\nu ^{(2)}$ is contained in a (natural) proper $G$-invariant integrable distribution.

\subsection {Homogeneous geodesics tangent to the nullity}\label {homogeo}
\
\smallskip

Let $M =G/H$ be  a presentation of a homogeneous Riemannian manifold, where $G$ is a connected Lie group which acts on $M$ by isometries. Fix $p\in M$ and let $N(p)$ be the leaf of nullity by  $p$.

If $X$ is a Killing field that is tangent to $N(p)$ at $p$, then $X_{\vert N(p)}$ must be always tangent to $N(p)$. This follows from the fact that  $X$ is projectable to the quotient of
$M$ by the integral manifolds of $\nu$ (see Section \ref{LTA}).

Let, for $p\in M$,
$$\gnul^p = \{g\in G: gN(p)=N(p)\}.$$
Then $\gnul^p$ is a Lie subgroup of $G$ which acts smoothly and transitively  on $N(p)$ (this action may be non-effective).

The Killing fields of $M$  induced by $\gnul^p$ are
those Killing fields  induced by $G$ that are tangent to $N(p)$ at $p$ (or equivalently, are always tangent to $N(p)$).  They form a Lie subalgebra of $\mathcal K(M)$ which we denote by
$\lienul^p \simeq \mathrm {Lie} (\gnul ^p)$.

Observe that the totally geodesic submanifold  $N(p)$ of $M$ is  extrinsically homogeneous and flat. Moreover, from  Section \ref
{homogeneous flat spaces}, it is globally flat.  Recall that for any Killing field $Z$ of $M$,
$\nabla Z$ is parallel along $N(p)$ (see Lemma \ref {ABC}). So, any  transvection
$X\in \mathcal K ^G(M)$ at $p$ must be also a transvection at all $q\in N(p)$.  If in addition it is tangent to $N(p)$, it is called an {\it extrinsic transvection} of $N(p)$.

\begin {lema}\label {SC} Assume that $M=G/H$ has non-trivial nullity distribution $\nu$ and does not split off, locally, a Euclidean factor. Then any leaf of nullity $N(p)$ is simply connected and so isometric to a Euclidean space. Moreover, the pullback  $i^*(TM)$ over the inclusion $i: N(p)\to M$ is globally flat.
\end {lema}

\begin {proof}
 Assume that $N(p)$ is not simply connected. We have, from Section \ref {homogeneous flat spaces}, that $N(p)$ has a non-trivial closed geodesic
$\gamma _v$ (in fact, this is a general fact about homogeneous spaces since any geodesic loop must be a closed geodesic). From Lemma \ref {ABC} any Killing field $X$ must be parallel along $\gamma _v$ and thus $\nabla _vX = 0$. Then, by Corollary  \ref {KostantdeRham}, $\mathbb R v$ extends locally to a parallel distribution and so $M$ locally splits off a line. A contradiction which proves the first assertion. The second assertion follows from the fact that $R_{u,v} =0$, if $u,v \in \nu _q$, and the first part.
\end {proof}

\

With the same arguments as in Section \ref {homogeneous flat spaces}, by considering isometries
$g_t \in \gnul ^p$ we have the following results

\begin {lema}\label {transvectionN(p)1} If there exists $X\in \lienul ^p$ such that
$X_{\vert N(p)}$ is not an intrinsic transvection of $N(p)$, then there exists
$0\neq Y \in \lienul ^p$ such that it is an extrinsic transvection at any point of $N(p)$ (i.e. $\mathfrak p ^p \cap \lienul ^p \neq \{0\}$).
\end {lema}

\begin {lema}\label {transvectionN(p)2} Let $\mathcal D$ be the parallel distribution of
$N(p)$ which is given by the directions of the extrinsic transvections and let
$\mathcal D^\perp$ be its  complementary perpendicular distribution on $N(p)$. Then, for any $q\in N(p)$, it holds:

{\it a}) If $v\in \mathcal D_q$,  then the geodesic $\gamma _v(t)$ is the integral curve of an extrinsic transvection $X \in \lienul ^p$ of $N(p)$.

{\it b}) If $v\in \mathcal D^\perp _q$, then the geodesic $\gamma _v (t)$ is the  integral curve of some $X\in \lienul ^p$.
\end {lema}

\

%%%%%%%%%%%%%%%%
%%%%%%%%%%%%%%%%
%%%%%%%%%%%%%%%%

Let $\hat \nu$ be the $G$-invariant distribution of $M$ defined by

\begin {equation} \label {hat}
\hat {\nu} _q  := \mathrm {span} \, \{ \nabla _wZ: Z\in \mathcal K ^G (M), w\in \nu _q\}.
\end {equation}

Then, from Lemma \ref {derived},
\begin {equation} \label {01nu} \nu ^{(1)} = \nu + \hat {\nu},
\end {equation}
where $\nu ^{(1)}$ is the osculating distribution.
Note that $\nu \subsetneq \nu ^{(1)}$, since $\nu$ is not a parallel distribution.

\begin {defi}\label {defiadapted} The distribution $\hat \nu$ will be called the {\it adapted distribution} of $\nu$.

\end {defi}

\begin {prop} [Existence of transvections] \label {existstransv}
Let $M=G/H$ be a  homogeneous Riemannian manifold, which does not split off (locally) a flat factor and with   a non-trivial nullity distribution $\nu$.
Then, for any $y\in \hat \nu _p$, there exists a transvection $Y\in \mathcal K ^G(M)$ such that $Y_p = y$.
\end {prop}
\begin {proof}

Let $Z\in \mathcal K ^G$ and let $v\in \nu _p$ that   either belongs  to $\mathcal D_p$ or  to $\mathcal D_p^\perp$ (we keep the notation of  Lemma \ref {transvectionN(p)2}). Then, by the above mentioned lemma, there exits $X\in
\lienul ^p$ such that $\gamma _v (t) = \phi _t (p)$, where $\phi _t$ is the flow associated to
$X$.

\color {black}

Let $Z\in \mathcal K^G(M)$ with $\nabla _vZ\neq 0$. Then, from Lemma \ref {ABC},
 $$ Z_{\gamma_v (t)} =
\tau_t (Z_p  + t\nabla _v Z),$$
where $\tau _t$ is the parallel transport along $\gamma (t)$, and
 $$ \nabla _{\gamma_v'(t)}(\nabla Z )=0,$$
  or equivalently
$$ (\nabla Z)_{\gamma_v (t)} = \tau _t ((\nabla Z)_p) =  \tau _t \circ (\nabla Z)_p\circ  \tau ^{-1}_t .$$

Let us  consider the family
$$Z^t: = (\phi _{-t})_*(Z)\ \ \ \ \ \ \ \ (t\in \mathbb R)$$of Killing fields induced by $G$.
Let us  compute their initial conditions at $p$. First  recall that from (\ref {etnabla}),
$\tau _t  ^{-1} \circ \mathrm d _p \phi _t = \mathrm e  ^{t(\nabla X)_p}$
 and so, since $(\mathrm d _p \phi _t )^{-1}= \mathrm d _{\gamma _v (t)}\phi _{-t}$,
 $$
 \mathrm d _{\gamma _v (t)}\phi _{-t}\circ \tau _t = \mathrm e  ^{-t(\nabla X)_p}$$
\begin{eqnarray*}Z^t_p &=& \mathrm {d}_{\gamma _v(t)}\phi _{-t} (Z_{\gamma (t)})
=  \mathrm {d}_{\gamma _v(t)}\phi _{-t} (\tau_t (Z_p  + t\nabla _v Z)) \\
&=&\mathrm e  ^{-t(\nabla X)_p}Z_p + t \mathrm e  ^{-t(\nabla X)_p}\nabla _vZ.
\end {eqnarray*}
\begin{eqnarray*}(\nabla Z^t)_p &=&
\mathrm {d}_{\gamma_v (t)}\phi _{-t} ((\nabla Z)_{\gamma (t)}) =
\mathrm {d}_{\gamma_v (t)}\phi _{-t} (\tau _t ((\nabla Z)_p)) \\
&=& \mathrm e  ^{-t(\nabla X)_p}((\nabla Z)_p) =
\mathrm e  ^{-t(\nabla X)_p}\circ (\nabla Z)_p \circ
\mathrm e  ^{t(\nabla X)_p}.
\end {eqnarray*}

Consider now the family $\frac 1 t Z^t$, $t\neq 0$. They are also Killing fields induced by $G$ and, by Remark \ref {accumulation}, we can choose  a sequence  of real numbers
$\{t_n\} \to +\infty$  such that $\mathrm e  ^{-t_n(\nabla X)_p}$ tends to the identity transformation of $T_pM$.
 Then $\frac {1\, } {t_n} Z^{t_n}$ converges
to a Killing field $Y$ with initial conditions (see Remark \ref{limitK})
$$(Y)^p = (\nabla _vZ, 0).$$

Since $\nu _p = \mathcal D_p \oplus \mathcal D^\perp _p $, and the sum of two transvections at $p$ is a transvection at $p$,  we finish the proof.
\end {proof}

\begin {nota}\label {corpropo} Since $\nu \subsetneq \nu ^{(1)} =
\nu + \hat {\nu}$, there must exist $Z\in \mathcal K ^G(M)$  and  $v\in \nu _p$ such that $v$ either belongs to $\mathcal D_p$ or to
$\mathcal  D^\perp _p$  and  $\nabla _vZ\notin \nu _p$.
Then, from Proposition \ref {existstransv},  there exists a non-trivial transvection
$Y$ at $p$ such that $Y_p = \nabla _vZ \notin \nu _p$.
\end {nota}

\
%begin  {nota}\label{otratransv}
%Let us keep the notation of
 %Corollary \ref {corpropo} and its proof.
 %Let $w=\nabla _vZ$ and assume that
 %$\nabla _wX\neq 0$.
 %Then the orbit $\mathrm {e}^{-t(\nabla X)_p}w$
 % in  $T_pM$ is non trivial. Choose $t_0$ such
 %that $w' = \mathrm {e}^{-t_0(\nabla X)_p}w$ is
 %linearly independent with $w$.
 %Following the proof of Proposition
 %\ref {existstransv}, and choosing a
 %$\{t_n\}_{n\in \mathbb N}\to \infty$ and such
 % that $\mathrm {e}^{-t_n(\nabla X)_p}$ tends to
 %$\mathrm {e}^{-t_0(\nabla X)_p}$ we construct
 % also a transvection $Y'$ at $p$ with $Y'_p=w'$
%\end{nota}

\begin {defi}\label {adapted} A transvection $Y$ at $p$ in the direction of
 $\nabla_v Z\in \hat \nu _p$, $Z \in \mathcal K ^G (M)$, $v\in \nu _p$,
will be called an {\it adapted transvection} to the vector $v$.
(See Proposition \ref {existstransv}).
\end {defi}

\

\subsection{Curvature of adapted transvections} \label{secmainteo}

\begin {prop}\label {R=0} Let $M=G/H$ be a homogeneous Riemannian manifold which does not split off a local flat factor.
Assume that $M$ has  a non-trivial nullity distribution  and let $Y$ be an adapted  transvection to the vector $v\in \nu _p$. Then, for any $U\in \mathcal K ^G(M)$  that  is bounded along $\gamma _v (t)$,  one has that
$$R_{U_p, Y_p}= 0.$$
\end {prop}

\begin{proof}
From Lemma \ref{SC} the integral manifold $N(p)$ of the nullity $\nu$ is complete and simply connected, and hence isometric to $\mathbb R ^m$.  Moreover, the pull-back
$i_* (TM)$  is globally flat.
Let  $\gamma_v (t)$ be a geodesic in $N(p)$, where $v\in T_pN(p) = \nu _p$. Write $v = v_1 + v_2$, where $v_1\in \mathcal D_p$,
$v_2\in \mathcal D^\perp _p$ (see Lemma \ref  {transvectionN(p)2}). Then the parallel transport $\tau _t$ along  $\gamma _v$, from $0$ to $t$, can be achieved as the composition of the parallel transport $ \bar \tau _t^t \circ \bar \tau _t$ where
$\bar \tau _t $ is the parallel transport along the  geodesic $\gamma _{v_2}(t)$, from $0$ to $t$,  and
$\bar \tau _t^t$ is the parallel transport along the geodesic $\gamma  ^t (s) = \gamma _{v_2}(t) + \gamma _{v_1}(s)$, from $s=0$ to $s=t$ (we identify $N(p)\simeq \mathbb R^m$).
Observe, that for any $t$,   there is a transvection $H^t$ such that its associated flow
$s\mapsto \psi^t_s$, satisfies
$\gamma ^t(s) = \psi^t_s(\gamma _{v_2} (t))$.

The geodesic $\gamma _{v_2}$ is homogeneous, i.e. $\gamma _{v_2} (t) = \phi  _t (p)$, where  $\phi _t$ is the flow associated to a Killing field $X\in \lienul ^p$, with $X_p= v_2$.

Then, from formula (\ref {etnabla}),
$$\bar \tau _t =\mathrm d _{p}\phi _t \circ \mathrm {e}^{-t(\nabla X)_p}, $$
 and
 $$\bar\tau ^t_t = \mathrm {d}_{\gamma _{v_2}(t)}\psi^t_t =
 \mathrm {d}_{\phi _t (p)}\psi^t_t $$
 Then

$$ \tau _t = \mathrm {d}_{\phi _t (p)}\psi^t_t
\circ \mathrm d _{p}\phi _t \circ \mathrm {e}^{-t(\nabla X)_p}. $$

If $g^t= \psi _t^t\circ \phi _t$, then $g^t$ is an isometry and

\begin {equation}\label {FFF}
 \tau _t = d_pg^t \circ \mathrm {e}^{-t(\nabla X)_p}.
 \end {equation}

Let $Z\in \mathcal K ^G(M)\simeq \mathfrak g $ be such that
$Z$ is not bounded along $\gamma _v (t)$, or equivalently, by Lemma \ref {ABC},
$w = \nabla _vZ \neq 0$.  Let $Y\in \mathcal K ^G(M)$ be a transvection adapted to $v$ whose initial conditions at $p$ are
$$(Y)^p = (w, 0).$$
% From Proposition \ref {existstransv}, such  a $Y$ does exist and it is a so-called transvection adapted to $v$.

 Let $U\in \mathcal K ^G (M)$ be bounded along $\gamma _v(t)$, or equivalently, $\nabla _v U =0$. We will determine the initial conditions of the bracket $[U,Z]$ at a point $\gamma _v (t)$, for an arbitrary $t$. More precisely, we are interested in the
 second component, $(\nabla [U,Z])_{\gamma _v (t)}$.
Recall, from Lemma \ref {ABC} (ii), that for any $\hat Z\in \mathcal K ^G(M)$, $\nabla \hat Z$ is parallel along $\gamma _v $.
From (\ref {initialbracket}), one has that

\begin {equation}  \label {UZ}
 (\nabla [U,Z])_{\gamma _v (t)} = R_{U_{\gamma _v (t)} , Z _{\gamma _v (t)}} - [(\nabla U)_{\gamma _v (t)}, (\nabla Z)_{\gamma _v (t)}].
\end {equation}

Since $\nabla U$ and $\nabla Z$ are parallel along $\gamma _v $, so is the bracket $[\nabla U, \nabla Z]$.
But $\nabla [U,Z]$ is parallel as well, and hence,  from (\ref {UZ}),
\begin {equation} \label {RUZ}
R_{U_{\gamma _v (t)} , Z _{\gamma _v (t)}}
\end {equation}
 must be   parallel along $\gamma _v$.

One has, from Lemma \ref{ABC}, that $U_{\gamma _v (t)} =
\tau _t (U_p)$ and that $Z_{\gamma _v (t)} =
\tau _t (Z_p) + t \tau _t w$, where $w = \nabla _vZ$. So, replacing in (\ref {RUZ}), one obtains that
\begin {equation}
R_{\tau _t (U_p), \tau _t (Z_p)} + t R_{\tau _t (U_p), \tau _t (w)}
\end {equation}must be parallel along $\gamma _v (t)$.
In particular, this expression must be bounded along $\gamma _v(t)$.

Since $M$ is homogeneous, both curvature operators $R_{\tau _t (U_p), \tau _t (Z_p)}$ and
$R_{\tau _t (U_p), \tau _t (w)}$ are bounded by
the supremum $$\mathrm{sup}\{\,\Vert R_{x,y}\Vert\,:\,x,y\in T_pM,\, \Vert x\Vert,\Vert y\Vert\leq C\}<\infty$$ for a suitable constant $C$.
This implies that $\Vert R_{\tau _t (U_p), \tau _t (w)}\Vert $ should tend to $0$ as $t$ tends to infinity.

Then, recalling (\ref {FFF}),

\begin {eqnarray} \label {RUZ4}
\Vert R_{\tau _t (U_p), \tau _t (w)}\Vert &=&
\Vert R_  { \mathrm d _{p}g^t (a(t)),
  \mathrm d _{p}g^t (b(t))}   \Vert
 \\ & =& \Vert  d _{p}g^t\circ R_  { a(t),
  b(t)}\circ \mathrm (d _{p}g^t)^{-1} \Vert \\
  &=& \Vert R_{a(t), b(t)}\Vert,
\end {eqnarray}
where
$a (t) = \mathrm {e}^{-t(\nabla X)_p}(U_p)$,
$b(t) = \mathrm {e}^{-t(\nabla X)_p}(w)$.

By Remark \ref {accumulation} one can take a sequence $\{t_n\}$ tending to infinity such that  $\mathrm {e}^{-t_n(\nabla X)_p}$
tends to the identity of $T_pM$. Then one concludes that
$\Vert R_{a(t_n), b(t_n)}\Vert $ tends to $ \Vert R_{U_p, w}\Vert$.
Since $\Vert R_{\tau _{t_n} (U_p), \tau _{t_n} (w)}\Vert =
\Vert R_{a(t_n), b(t_n)}\Vert $ must tend to $0$ as $t\to \infty$,
we conclude that  $R_{U_p, w}=0$. \end{proof}

\begin {nota}\label {generalizedR=0} In the proof of  Proposition \ref {R=0}  we only used that the projection  of $U$ to $\nu ^\perp$ is bounded. Then: {\it $ R_{U_p, Y_p}= 0$ if  $Y$ is an adapted transvection at $p$ and  $U$ belongs to the bounded algebra $\mathfrak u ^p$}
\noindent (see Definition \ref {BA} and its preceding paragraph).
\end {nota}

\begin {lema}\label {mainL}
Let $M=G/H$ be a homogeneous Riemannian manifold which does not split off a local flat factor and with a non-trivial nullity distribution.   Let   $0\neq Y$ be a
transvection adapted to $v\in \nu _p$.
\color {black} Then
\begin {enumerate}
\item If $Z$ is any Killing field of $M$, then $[Y,Z]$ is bounded along $\gamma _v$ (or equivalently, $\nabla _v[Y,Z] = 0$).
\item If $U$ is any  Killing field of $M$ which is bounded along $\gamma _v$, then
$[Y, U]$ is a transvection at $p$ (i.e., $(\nabla[Y, U])_p= 0$).
\item If $\bar Y$ is any transvection at $p$, then
$[Y, \bar Y] =0$.
\item
$[Y,[Y,[Y, \mathcal K (M)]]]= 0$, or equivalently,
identifying Killing fields with elements of the isometry algebra,
$\mathrm {ad}^3_Y = 0$, in the Lie algebra of the full isometry group of $M$.
 \end{enumerate}
\end {lema}

\begin {proof}
Let $Z\in \mathcal K (M)$. Then,   from
(\ref{initialbracket}), $(\nabla [Y,Z])_p = R_{Y_p, Z_p}$.
So, $$\nabla _v [Y,Z] = R_{Y_p, Z_p}v=0 ,$$ since $v$ is in the nullity.
Then $[Y,Z]$ is bounded along $\gamma _v$. This proves (i).\\

To see (ii) observe that $(\nabla [Y,U])_p = R_{Y_p, U_p} = 0$ by Proposition \ref {R=0}. Thus, $[Y,U]$ is a transvection at $p$.\\

Let now $\bar Y$ be a transvection at $p$. In particular $\nabla _v\bar Y =0$ and so $\bar Y$ is bounded along $\gamma _v$. Then part (ii) applies and $(\nabla [Y, \bar Y])_p = 0$. Since $Y, \bar Y$ are both transvections at $p$, by (\ref {initialbracket}), $ [Y, \bar Y]_p =0$.
Then $([Y, \bar Y])^p =(0,0)$ and so $[Y, \bar Y]$ vanishes identically on $M$. This proves (iii).\\

Finally if $Z\in \mathcal K (M)$ is arbitrary, then $[Y,Z]$ is bounded
by part (i). Then applying (ii) $[Y,[Y,Z]]$ is a transvection at $p$.
Then, by part (iii), $$[Y,[Y,[Y,Z]] = 0.$$
\end {proof}

We improve part (iv) of Lemma \ref {mainL}. Namely,

\begin {teo}\label {ST}
Let $M=G/H$ be a homogeneous Riemannian manifold which does not split off a local flat factor and with a non-trivial nullity distribution.   Let   $0\neq Y$ be a transvection adapted to $v\in \nu _p$. Then
$[Y,[Y, \mathcal K (M)]]=0$, or equivalently, identifying Killing fields with elements of the isometry algebra,
$\mathrm {ad}^2_Y=0$, in the Lie algebra of the full isometry group of $M$. Moreover, $[Y, \mathcal K ^G(M)]\neq 0$.
\end{teo}

\begin {proof}
We will not regard, as before, Killing fields along $\gamma _v$, but along the geodesic $\beta (t) = \phi _t (p)$, where $\phi _t$ is the flow associated to $Y$.
Since $Y$ is a transvection at $p$, $d_p\phi _t$ coincides with the parallel transport along $\beta (t)$. Then if $\psi$ is any field in $M$,
$[Y,\psi]_{\beta (t)}$ is the covariant derivative $\frac {\, \mathrm D}
{\mathrm d t} \psi _{\beta (t)}$. Let us apply this for $\psi = Z\in \mathcal
K (M)$. Keep in mind that $Z_{\beta (t)}$ is a Jacobi field along $\beta$. So, from
Lemma \ref {mainL}, (iv)
\begin {equation}\label {D3}
\frac {\, \mathrm D^3}
{\mathrm d t^3} Z _{\beta (t)} =0.
\end {equation}

Now in general the curvature tensor $R$ is invariant under isometries and
$d_p\phi _t$ coincides with the parallel transport along $\beta (t)$.
This implies that the Jacobi operator $R_{\, \cdot , \beta' (t)} \beta' (t)$
diagonalizes in a parallel basis with constant distinct  eigenvalues $\lambda _0 =0, \lambda _1, \cdots , \lambda _r$ (as in symmetric spaces (see \cite {BOR})).

Let  $V^0 (t), V^1 (t), \cdots , V^r(t)$ be the  eigenspaces of the Jacobi operator $R_{\, \cdot , \beta' (t)} \beta' (t)$ associated to
$0, \lambda _1, \cdots , \lambda _r$, respectively. Any of such subspaces must be  parallel along $\beta (t)$.
Then the orthogonal projection $Z^i(t)$ of $Z_{\beta (t)}$ to
 $V^i(t)$ is of one of the following  types, according with the sign of $\lambda _i$.

 (a) $Z^0(t) = a(t) + t b(t)$,
 where $a(t), b(t)\in V^0(t)$ are parallel fields along $\beta (t)$.

 (b) If $\lambda _i >0$,
 $Z^i(t) = \cos (\sqrt {\lambda_i} t) a(t) + \sin  (\sqrt {\lambda _ i }t) b(t)$, where  $a(t), b(t)\in V^i(t)$ are parallel fields along $\beta (t)$.

  (c) If $\lambda _i <0$, $Z^i(t) = \cosh (\sqrt {-\lambda_i} t) a(t) + \sinh (\sqrt {-\lambda_i} t) b(t)$, where  $a(t), b(t)\in V^i(t)$ are parallel fields along $\beta (t)$.\\

But (\ref {D3}) implies that $Z^i(t) = 0$, for $i= 1, \dots , r$ hence $Z _{\beta (t)} = Z^0(t)$ and so
\begin{equation}\label {D4}
\frac {\, \mathrm D^2}
{\mathrm d t^2} Z _{\beta (t)} =0.
\end {equation}

%Being $M$ homogeneous, there is a   Killing fields in any direction, we conclude that the Jacobi operator has only one eigenvalue $\lambda _0 = 0$.
%
%The Jacobi equation implies that
%\begin{equation}\label {D4}
%\frac {\, \mathrm D^2}
%{\mathrm d t^2} Z _{\beta (t)} =0
%\end {equation}

Then $[Y,[Y, Z]]_{\beta (t)}\equiv 0$. In particular,
$[Y,[Y, Z]]_p = 0$. But,  by Lemma \ref {mainL} (i) and (ii),
$(\nabla[Y, [Y,Z]])_p = 0$. Then  $[Y,[Y, Z]] =0$.
This proves  that $$[Y,[Y, \mathcal K (M)]]=0 .$$

It only remains to show that $[Y, \mathcal K ^G(M)]\neq 0$.
Assume, on the contrary, that  $[Y, Z]= 0$, for all $Z\in \mathcal K^G(M)$.
Since $Y$ is a transvection at $\beta (t)$, for all $t$,
the covariant derivative along $\beta (t)$  of
$Z_{\beta (t)}$, as we have seen, coincides with $[Y, Z]_{\beta (t)} =0$ as assumed.
Then $Z_{\beta (t)}$ is parallel along $\beta (t)$ hence
$\nabla _{\beta'(0)}Z = \nabla _{Y_p}Z = 0$. Since $Z$ is arbitrary in $\mathcal K ^G(M)$ we conclude, from Corollary \ref {KostantdeRham}, that $M$ splits locally the direction of $Y_p$. A contradiction.
\end {proof}

In the proof of the above theorem it was shown that the adapted transvection $Y$ has null Jacobi operator along $\beta (t)$ or equivalently at $p$. Indeed being $M$ homogeneous, there is a Killing field in any direction, and we conclude that the Jacobi operator has only one eigenvalue $\lambda _0 = 0$. From Remark \ref {corpropo} we may assume that
$Y_p \notin \nu _p$. Then we have the following result that will be very useful for finding irreducible homogeneous Riemannian manifolds with non-trivial nullity distribution.

\begin {cor}\label {STcor} Let $M=G/H$ be a homogeneous Riemannian manifold which does not split off a local flat factor and with a non-trivial nullity distribution $\nu$.   Then any adapted transvection $Y$ to $v\in \nu _p$ has a null Jacobi operator $R_{\, \cdot , Y_p}Y_p$.
Moreover, there exists such an adapted  transvection  $Y$ with $Y_p\notin \nu _p$.
\end {cor}

Any transvection at $p$ belongs, by definition, to the Cartan
subspace $\mathfrak p ^p$ at $p$ (see (\ref {Cs})).
Those with trivial Jacobi operator must lie in the abelian part $\mathfrak
p ^p_0$ of the Cartan subspace. Namely,

\begin {cor}\label {ABE} Let $M=G/H$ be a homogeneous Riemannian manifold which does not split off a local flat factor and with a non-trivial nullity distribution. Then any transvection $Y$ adapted to  $v\in \nu _p$ belongs to the abelian part
$\mathfrak p^p_0$ of the Cartan subspace at $p$. In particular, the distribution of symmetry $\mathfrak s$ of $M$  is non-trivial and so the index of symmetry of $M$ is positive.
\end {cor}

\subsection {Transvections with null Jacobi operator}\label {NuJaOp} \
\smallskip

We have the following result for homogeneous spaces that have transvections with null Jacobi operator (not depending on the existence of a non-trivial nullity).

\begin {prop} \label {JacNul} Let $M=G/H$ be a homogeneous Riemannian manifold where $G = I(M)^o$ and such that $M$ does not split off, locally, a flat factor. Let $\mathfrak a_0^p$ be the set of transvections $Y$ at $p\in M$ with null Jacobi operator $R_{\cdot , Y_p}Y_p$. Then
\begin {itemize}
\item [(i)] If $0\neq Y\in \mathfrak a_0^p$, then $[Y,[Y, \mathcal K (M)]] = 0$ and
$[Y, \mathcal K^G(M)]\neq 0$.

\item [(ii)] The set $\mathfrak a_0^p$ is a vector space. Moreover, it is an abelian Lie algebra and the distribution $q \to  \mathfrak a_0^q.q$ is autoparallel and flat.
\end {itemize}
\end {prop}

\begin {proof}

Let $Y\in \mathfrak a_0^p$. Let $\phi _t$ the flow associated to $Y$. Then
$\gamma (t)  =\phi _t(p)$ is the geodesic with initial condition $v=Y_p$. Let $Z\in \mathcal K (M)$. Then, from the Jacobi equation,
$Z_{\gamma (t)} = u(t)  + t w(t)$, where $u(t), w(t)$ are the parallel transports along $\gamma (t)$ of $Z_p$ and $\nabla _vZ$, respectively. Since
$[Y,Z] = \frac {\mathrm \, d}{\mathrm d t}_{\vert 0}
(\phi _{-t})_*(Z)$, by formula
(\ref {etnabla}), we obtain that
$$[Y, Z]_{\gamma (t)} =\nabla _{\gamma ^{\prime}(t)}Z_{\gamma (t)}=  w(t), $$
i.e., $[Y, Z]_{\gamma (t)}$ is parallel along $\gamma (t)$.
So,
$$ [Y,[Y, Z]]_{\gamma (t)}=0. $$
 Then $U =[Y,[Y, Z]]$ belongs to the double isotropy algebra at $v$. In particular
$U$ belongs to the isotropy algebra $\mathcal{K}(M)_p \subset \mathfrak {so}(T_pM)$
(via the isotropy representation at $p$). Moreover, if $\psi _s$ is the flow associated to $U$, $\mathrm d _p \psi _s (v)= v$, for all $s$. Or, equivalently,
$$\psi _s (\gamma (t)) = \gamma (t).$$
Since isometries map transvections into transvections,  from the above equalities we obtain that $(\psi _s)_* (Y) = Y$ and therefore, differentiating the flow,  $-[U,Y]= [Y,U]=0$.

%The same argument shows that $[W,Y] =0$,
%for all $W$ in the double isotropy
%$(\mathfrak g _p)_v\subset \mathfrak g _p$
%at $v$.

Let $B$ be the Killing form of $\mathcal K(M)$. Then $B_{\vert \mathcal{K}(M)_p\times \mathcal{K}(M)_p} $
is negative definite (see the beginning of the proof of Lemma \ref {EFF}).
$$B(U,U) = B([Y,[Y, Z]],U) = - B([Y, Z],[Y,U]) = 0.$$
Then $U=0$, since $U\in \mathcal{K}(M)_p$. This proves the first assertion of (i).

If $[Y, \mathcal K^G(M)] =  0$, then $g_*(Y) = Y$, for all $g\in G$. Then, since $G$ acts transitively on $M$,  $Y$ is a transvection at any point and so a parallel field.
Then $M$ locally splits off a line. A contradiction. This finishes the proof of (i).

Since $\mathfrak a_0^p$ is invariant by scalar multiplications, one has to show that if $X, Y\in \mathfrak a_0^p$ then $Z= X+ Y \in \mathfrak a_0^p$
(observe that the Jacobi null condition is not a linear condition). In the notation of
Section \ref {index}, we have that $\mathfrak a_0^p \subset \mathfrak p ^p_0$, the abelian part of the Cartan subspace at $p$. Observe that $[X, Y]=0$, since this field is zero when restricted to the symmetry leaf $L(p)$ (see Lemma \ref {EFF}).
Then $[X, [Y, \cdot]] = [Y, [X, \cdot]]$. Taking this and the fact that
$[X, [X, \cdot]] = 0 = [Y, [Y, \cdot]] $ into account, we obtain that
$$[Z, [Z, [Z, \cdot]] = 0.$$
The same argument used in the proof of Theorem \ref {ST} shows that the Jacobi operator $R_{\cdot , Z_p}Z_p$ is zero (see also the paragraph just after this theorem). Then $Z\in \mathfrak a^p_0$. So $\mathfrak a_0^p$ is a vector subspace of the abelian Lie algebra $\mathfrak p _0 ^p$. If $L(p)$ is the leaf of symmetry by $p$, then its Euclidean factor is given by $A.p$, where $A$ is the abelian connected Lie group of isometries with Lie algebra $\mathfrak p _0 ^p$. If $A'$ is the connected Lie group associated to $\mathfrak a^p_0$, then $A'.p$ is totally geodesic in $A.p$. Since $A.p$ is totally geodesic in $M$, then $A'.p$ is totally geodesic in $M$ and it is an integral manifold of  the distribution $q \to  \mathfrak a_0^q.q$. This concludes the proof.
\end {proof}

Observe, from Proposition \ref {JacNul}  and Theorem \ref {ST}, that in any direction of $\hat \nu_p$ there exists a transvection
$X$ with $\mathrm {ad}^2_X = 0$ and $\mathrm {ad}_X \neq 0$.

\vspace{.2cm}

\color {black}
\subsection{Applications to semisimple and nilpotent homogeneous spaces}$\ $

\begin{cor}\label{corsemisimple}
Let $M=G/H$ be a simply connected homogeneous Riemannian manifold without a Euclidean de Rham factor.
Assume that Lie algebra $\mathfrak g$ of $G$ is reductive.
Then $M$ has a trivial nullity.
\end{cor}

\begin{proof}
Assume that the nullity $\nu$ is non-trivial and let $p= [e] \in M$.
We will regard $\mathfrak g$ as  the Lie algebra
$\mathcal K^G$ of Killing fields induced by $G$.
Then, by Proposition  \ref {existstransv}  and Theorem \ref {ST},
there exists an adapted transvection
$0\neq Y\in \mathfrak g$ of $M$ at $p$, such that
$\mathrm {ad}^2 _Y= 0$ and $\mathrm {ad} _Y\neq  0$. Let us decompose
$\mathfrak g = \mathfrak a \oplus \mathfrak k \oplus \mathfrak i$
into the direct sum the ideals which are abelian, of the compact type, and of the  noncompact type, respectively.

Let us write $Y = Y^0+
Y^1 + Y^2$, where $Y^0\in \mathfrak a$,  $Y^1\in \mathfrak k$, and
 $Y^2\in \mathfrak i$. Since $\mathrm {ad}^2 _Y= 0$, then $\mathrm {ad}^2 _{Y^1}=0$. Hence $Y^1=0$, since $\mathfrak k$ is of the compact type.

 On the one hand, since $\langle \nabla _{\mathfrak g}Y,\mathfrak a\rangle_p = 0$, we obtain, from equation (\ref{fundamentalequation}), that
 \begin {equation}\label {perAb}
 \langle [Y, \mathfrak g], \mathfrak a\rangle _p =0.
 \end {equation}

On the other hand,  since $\mathrm {ad} _Y\neq  0$,  $Y^2\neq 0$. Hence $Y^2$ is a 2-step nilpotent element of $\mathfrak i$.
Then, by the Jacobson-Morozov theorem \cite{OVG}, $Y^2$ belongs to a $\mathfrak {sl}(2)$-triple,  in $\mathfrak i$, $\{Y^2, Z,W\}$ such that
$[W,Y^2]=2Y^2$, $[W,Z]=-2Z$ and $[Y^2,Z]=W$.

Since $Y$ is a transvection at $p$, by equation
(\ref{fundamentalequation}), we have that
\begin {eqnarray*} 0&=&\langle \nabla_Y Y,W\rangle_p=\langle[Y,W],Y\rangle_p
= \langle[Y,W],Y ^2\rangle_p \text {\small \ \ by (\ref {perAb})}  \\
&=& \langle[Y^2,W],Y^2\rangle_p=
-2\Vert Y^2_p\Vert^2.
\end {eqnarray*}
Then $Y^2_p=0$ and so $Y^2$ belongs to the isotropy algebra
$\mathfrak h = \mathrm {Lie} (H)$. Let $(\, , \, )$ be an $\mathrm {Ad}(H)$-invariant inner product in $\mathfrak g$. Then $\mathrm {ad} _{Y^2}$ is skew-symmetric and so
$\mathrm {ad} _{Y^2}=0$, since $\mathrm {ad}^2 _{Y^2}=0$.
Then $Y^2\in \mathfrak a$. A contradiction.
\end {proof}

\bigskip

The same argument of the above corollary shows that there are no transvections of order $2$ in a homogeneous Riemannian manifold $M=G/H$ with $\mathfrak g$ reductive.

\

\color {black}

\begin{cor}\label{cornilpotent}
Let $M = G$ be a 2-step nilpotent Lie group with a left invariant metric
which does not split off a local flat factor.
Then the nullity distribution of $M$ is trivial.
\end{cor}

\begin{proof} \label {2-step}
We proceed by contradiction assuming that the nullity distribution $\nu$
is non-trivial.
Let $Y$ be the transvection at $p \in M$ given by Theorem  \ref{ST} such that $[Y, \mathfrak g]\neq {0}$ i.e. $Y$ does not belong to the center $\mathfrak c$ of $\mathfrak g$. From (\ref
{fundamentalequation}),
$$0 =2\langle \nabla _{\mathfrak c} Y , \mathfrak g \rangle_p =
\langle [Y, \mathfrak g], \mathfrak c\rangle _p \, \, .$$

Since $[Y, \mathfrak g] \subset \mathfrak c$ due to the fact that $G$ is 2-step nilpotent we get that $[Y,\mathfrak g]_p = 0$
hence $[Y, \mathfrak g] = {0}$. A contradiction.
\end{proof}
%
%there is a non-trivial transvection at any $p\in M$ by Proposition \ref {existstransv}. Moreover, $[Y, \mathfrak g]\neq {0}$ and so $Y$ does not lie in  $\mathfrak c$.

\begin{nota} The above corollary also follows from well-known facts about the Ricci tensor and the de Rham factor of a 2-step nilpotent Lie group with a left invariant metric \cite[Proposition (2.5) and Proposition (2.7)]{Eb94} .
\end{nota}

\subsection {The nullity is not a homogeneous foliation} \

We finish this section by showing that the nullity foliation is far from being a homogeneous foliation (i.e. given by the orbits of  an isometric group action).

Observe, from the affine Killing equation (\ref {affK}), that a Killing field $X$ lies in the nullity
distribution if and only if the Nomizu tensor $\nabla X$ is parallel. Below we show that a Killing field $X\neq 0$ that lies in the nullity  must be  tangent to the (local) Euclidean de Rham factor. We were unable to find  this result in the literature.

\begin {prop}\label {KillTanNu} Let $M$ be a homogeneous Riemannian manifold without Euclidean (local) de Rham factor.  Assume that the nullity distribution $\nu $ of $M$ is non-trivial. Then there exists no  Killing field of $M$, not identically zero, such that $X$ is always tangent to $\nu$.
\end {prop}

\begin {proof}
We may assume that $M$ is simply connected. Let
$M = M_1 \times \cdots \times M_r$ be the de Rham decomposition of $M$, where
$M_1, \cdots , M_r$ are irreducible Riemannian manifolds. The nullity $\nu$  of $M$ is the direct sum of the nullities  $\nu _i$  of $M_i$, $i=1, \cdots , r$. Let $X\in \mathcal K(M)$ be always tangent to $\nu$. The projection  $X_i$ of $X$ to any given factor $M_i$ belongs to
$\mathcal K (M_i)$. If $M_i$ has a trivial nullity, then $X_i=0$.  So we may  assume  that
 $\nu _i \neq 0$.  Assume that $X_i\neq 0$. From formula (\ref {affK}) we have that $\nabla X_i$ is a parallel
 skew-symmetric $(1,1)$-tensor of $M_i$. If $\ker (\nabla X_i) =TM$, then $X_i$ is a parallel field and so $M_i$ splits off a line. A contradiction.   Let $p\in M$. Then $(\nabla X_i)_p$ has a complex eigenvalue  $\lambda\notin \mathbb R$. Let $\{0\}\neq \mathbb V\subset T_pM$ be the $(\nabla X_i)_p$-invariant  subspace
 associated to $\lambda$ and $\bar {\lambda}$. Since $\nabla X_i$ is a parallel tensor, then $\mathbb V$ extends to a parallel distribution of $M_i$ and since $M_i$ is irreducible, then $\mathbb V = T_pM$. In this case, eventually by rescaling $X$, we have that $J=\nabla X_i$ is a K\"ahler structure on $M_i$ and so $M_i$ is a K\"ahler (homogeneous) manifold. Then the field $\xi = JX_i$ lies in the nullity distribution $\nu _i$.
 Moreover, as it is standard to check, $\nabla ^2\xi = 0$. So
 $\xi$ satisfies the affine Killing equation
 $\nabla_{u,v} ^2 \xi = R_{u,\xi}v$ (see \ref {affK}). If $\phi _t$ is the flow associated to $\xi$, for any given $t$, $\phi _t$ is a homothetic transformation of $M$  associated to the constant $\mathrm {e}^{ta} $, where $A = a Id$ is the symmetric part of
 $\nabla \xi$ (cf.  \cite {KN}, Lemma 1, pg. 242). In fact, this symmetric part  $A$, from the affine Killing equation, must be parallel.  Moreover, since  $M_i$ is irreducible, $A$ has only one (constant) eigenvalue.
 In our particular case $\nabla \xi = J \nabla X_i= -Id $, and so $a=-1$.
 But,  in a homogeneous non-flat irreducible space,  any homothetic transformation is an isometry. This is a general fact for a  complete Riemannian manifolds (see \cite {KN}, Theorem 3.6, pg. 242). For the sake of self-completeness, we will show this in our homogeneous context. In fact,
 $\mathrm {d}_p \phi _t : T_pM \to T_{\phi _t(p)}M$ is a homothetic map. Namely,
 $$\langle \mathrm {d}_p \phi _t (v),  \mathrm {d}_p \phi _t (v)\rangle =
 \mathrm {e}^{-2t}\langle v, w\rangle.$$
   One has, since $\phi _t$ preserves the Levi-Civita connection,
   that  $\mathrm {d}_p \phi _t$ maps  the curvature tensor $R^p$ of $M_i$ at $p$ into the curvature tensor $R^{\phi _t(p)}$ at
 $\phi _t(p)$.  Then, by a standard calculation,
 $$\mathrm {e}^{-2t}\Vert R^{\phi _t(p)}\Vert = \Vert R^{p}\Vert.$$
 This is a contradiction, since $M_i$ is homogeneous and non-flat.  Then $X_i =0$, for any $i=1, \cdots , r$. Then $X=0$.
\end {proof}

The same proof works assuming $M$ to be complete not necessarily homogeneous.

\vspace {.2cm}

\section {Symmetry and nullity}\label{symnul}

Let $M=G/H$ be a homogeneous locally irreducible Riemannian manifold with a non-trivial distribution of  symmetry $\mathfrak s$. Recall that $\mathfrak s$ is not contained in the nullity distribution $\nu$, see Remark \ref {corpropo}.  Since both distributions $\nu$ and $\mathfrak s$ are $G$-invariant their sum
\begin {equation} \label {proper2}
\tilde \nu = \nu + \mathfrak s
\end {equation} has constant rank, and hence $\tilde \nu$ is a distribution on $M$.
Observe that the above sum could be non-direct.

\begin{lema}\label {proper}
The distribution $\tilde \nu$ is autoparallel.
Moreover, if $\tilde N (p)$ is an integral manifold of $\tilde \nu$ then
the restrictions $\mathfrak s_{\vert \tilde N (p)}$, $\nu_{\vert \tilde N (p)}$ are  parallel distributions of
$\tilde N (p)$.
\end{lema}

\begin{proof}

Let $Y\in \mathfrak p^p$, the Cartan subspace at $p$ (see (\ref {Cs})), and let $c(t)$ be a curve contained in the leaf of nullity $N(p)$ joining $p$ and an arbitrary point $q \in N(p)$.
From the affine Killing equation (\ref {affK}) one has that $\nabla Y$ is parallel along $c(t)$. This implies that
$(\nabla Y)_q = 0 $ for all $q\in N(p)$ hence $Y_q\in \mathfrak s _q$.
Since $p$ is arbitrary, we get
$$\nabla _{\nu}\mathfrak s \subset \mathfrak s .$$

Let $\phi _t$ be the flow associated to $Y$.
Since $\nu$ is $G$-invariant and, by equation (\ref {etnabla}),   $d_p\phi _t$ gives the parallel transport along (the geodesic) $\phi _t (p)$, we must have that
$\nu$ is parallel along the leaf of  symmetry $L(p)$ at $p$.
Since $p$ is arbitrary we conclude that
$$\nabla _{\mathfrak s} \nu \subset \nu  .$$
Then, since $\nu $ and $\mathfrak s$ are both autoparallel, we conclude that $\tilde \nu$ is autoparallel.
\end{proof}

Let now $\mathfrak s ^0 $ be the flat part of the distribution of symmetry (see equation (\ref{flatfactor})) and consider the distribution
   \begin {equation}\label {proper0}
 \tilde \nu ^0 = \mathfrak s ^0 + \nu
  \end {equation}
which is not in general a direct sum.

\

We have the following lemma.

\begin{lema} The nullity distribution is  properly  contained
in the $I(M)$-invariant distribution $ \tilde \nu ^0 $,
which is autoparallel and flat.
\end{lema}

\begin{proof}
First observe that, from Corollary \ref{ABE}, there is a transvection $Y\in\mathfrak{s}^0$ which does not lie in $\nu$. So $\nu$ is properly contained in $\tilde{\nu}^0$. Since both $\mathfrak{s}^0$ and $\nu$ are $I(M)$-invariant, so is $\tilde{\nu}^0$.

By the above Lemma we have that locally $\tilde N (p) = L(p) \times W$ as Riemannian product, where $W$ is a Riemannian manifold.
Now $\mathfrak s ^0 $ is the flat parallel distribution tangent to the whole flat de Rham factor of any leaf of symmetry $L(q)\subset \tilde N (q)$. So we conclude that the restriction $\mathfrak s ^0 _{\vert \tilde N (p)}$  is parallel.

% In fact, the parallel transport in $\tilde N (q)$ along a curve from $x$ to $y$, induces a local isometry between  $L(x)$ and $L(y)$ which maps the Euclidean de Rham distribution $\mathfrak s ^0_{\vert L(x)} $ of  $L(x)$ in the similar object in $L(y)$.

 Then
 $$\tilde \nu ^0_{\vert \tilde N (p)} = \mathfrak s ^0 _{\vert \tilde N (p)}+ \nu_{\vert \tilde N (p)}$$
 is a parallel distribution of $\tilde N(p)$. This implies that $\tilde \nu ^0 $ is an autoparallel distribution of $M$. Moreover, it must be flat, since
$ \mathfrak s ^0 $ and $\nu$ are parallel and flat distributions of
$ \tilde N (p)$.
\end{proof}

%{\color{red}
%\begin{nota}
%This result can be also obtained as a consequence of Proposition 4.23, \cite{BVK}.
%\end{nota}}

\section {The osculating distributions and the isotropy }\label{section5}

Here we give the details of the proofs of part (1) and (3) of Theorem A. Without loss of generality, keeping the notation of this theorem, we may assume that the presentation group $G$ is connected.

We already showed that
\begin {equation}\label {123} \nu ^{(1)} = \nu + \hat {\nu}
\end {equation}
and that $ \nu \subsetneq \nu ^{(1)} = \nu + \hat {\nu}$, see equation (\ref{01nu}).
Moreover, $\hat {\nu}\subset  \mathfrak s^0$ by Corollary \ref{ABE}. By using that $\hat {\nu}$  is $G$-invariant and making the same arguments as in Section \ref {symnul} one has that $\nu ^{(1)}$
is an autoparallel and flat, and so a proper, distribution of $M$.
One can also prove this fact by using that $R_{X,Y} = 0$, if $X, Y$ are transvections  that belong to $\mathfrak p _0^p$ (and $\hat {\nu}_p \subset \mathfrak p _0^p.\, p = \mathfrak s _0^p$).
\vspace{.25cm}

The inclusion $\nu ^{(1)} \subsetneq \nu ^{(2)}$ is proper. In fact, if $\nu^{(1)}=\nu^{(2)}$, $\nu ^{(1)}$ would be a parallel distribution and, since it is flat, $M$ would have a Euclidean de Rham factor, which contradicts the assumptions of Theorem A.

The osculating distribution $\nu ^{(2)}$ is $G$-invariant and by (\ref {123}) one has that
 $$\nu ^{(2)} = \nu ^{(1)} + \hat {\nu}^{(1)}.$$

Then, from Lemma \ref{derived}, and taking into account that $\hat {\nu}\subset
\nu ^{(1)}$, it follows that

\begin {equation} \label {nu2span}
\nu ^{(2)}_p =  \nu ^{(1)}_p + \mathrm {span} \, \{\nabla _v Z: Z\in \mathcal K ^G (M) , v\in \hat {\nu} _p\}.
 \end {equation}

 But, if $X$ is a transvection with $X_p =v\in \hat {\nu}_p$, $\nabla_vZ = [X,Z]_p+\nabla _{Z_p}X=[X,Z]_p$. Then

 \begin {equation} \label {biss}
\nu ^{(2)}_p =  \nu ^{(1)} +\mathrm {span} \,  \{[X,Z]_p: Z\in \mathcal K ^G (M), X \text { is a transvection at } p \text { with } X_p \in
 \hat {\nu} _p\}
 \end {equation}

 Let $X$ be a transvection at $p$, with $X_p\in \hat {\nu}_p$, and let
 $Z\in \mathcal K ^G (M)$ be arbitrary. Then, from   (\ref {initialbracket}), the initial conditions at $p$ of $[X,Z]$ are

 \begin {equation} \label {br}
 ([X,Z])^p = ([X,Z]_p , (\nabla [X,Z])_p) = (\nabla _{X_p}Z , R_{X_p,Z_p}).
 \end {equation}

 Observe that
 \begin {equation}\label {uuuu}
 \nabla_{\nu _p} [X,Z] = R_{X_p,Z_p}\nu _p = \{0\}.
 \end {equation}

 % \nu ^{(1)} + \mathrm  {span }
 %\{\nabla _vX: \, X \text {is a field that lie in }
% \hat {\nu} \text { and\ } v \in T_pM\}

\

 For $p\in M$, define the subspace $\mathfrak u^p$ of $\mathcal{K}^G(M)$ as $$\mathfrak u^p=\{ U\in \mathcal{K}^G(M)\,:\, \nabla_{\nu_p}U\subset \nu_p \}.$$
 From Lemma \ref{ABC} it follows that for any $q\in N(p)$, $\mathfrak u^p=\mathfrak u^q$. Moreover, from formula (\ref {curvatureformula}) it follows that $\mathfrak u^p $ is a Lie subalgebra of $\mathcal{K}^G(M)$. Observe that  $U\in \mathfrak u^p$ if and only if the normal component  to $\nu$ of $U_{\vert N(p)}$ is parallel with respect to the normal connection of $N(p)$.  Notice that such normal component is parallel with respect the normal connection if and only if it is parallel as a section of the pullback bundle $i_*(TM)$ endowed with the Levi-Civita connection, since $TN$ is a parallel sub-bundle of the pullback $i_*(TM)$, (cf. Lemma \ref{SC}).

  If $g$ is an isometry, then $\mathfrak u ^{g(p)}=
 g_*(\mathfrak u ^p)$.
 \

 \begin {defi}\label {BA} The Lie algebra $\mathfrak u ^p$ is called the {\it  bounded algebra at $p$.} The $G$-invariant distribution $\mathcal U$ defined by
 $\mathcal U_p = \mathfrak u^p .\, p$ is called the {\it  bounded distribution}.
 \end {defi}

 \

\noindent
{\bf Observation.} The Lie algebra  $\mathfrak u ^p$ contains:
\begin {itemize}
\item  Any transvection at $p$.

\item  Any element in the linear span of
$$\{[X,Z]: Z\in \mathcal K ^G (M), X \text { a transvection at } p \text { with } X_p \in  \hat {\nu} _p\}.$$

\item  Any Killing field which is tangent to $N(p)$ at $p$ (and so always tangent to $N(p)$). In particular, any Killing field in the isotropy  algebra at $p$.

\item The bounded distribution does not depend on the presentation group $G$ of $M$. This follows from the fact the bounded algebra contains the isotropy algebra.

\end {itemize}
\rm

\

From the above properties we have that
\begin {equation} \label {332}
\nu ^{(2)}_p \subset \mathfrak u ^p.\, p.
\end {equation}

% Let $\mathcal U$ be the $G$-invariant distribution of $M$ defined by
%
% $\mathcal U _p = \mathfrak u ^p .\, p$.

Let $\bar G^p \subset G$ be the Lie subgroup associated to the Lie subalgebra
$\mathfrak u ^p \subset \mathfrak g =
\mathrm {Lie}(G) \simeq \mathcal K^G(M)$.

Then, from the previous observation,  $\bar G^p$ contains the isotropy subgroup $G_p$. In fact, since $M$ is simply connected, and $G$ is assumed to be connected, then $G_p$ is connected.
Note  that $\bar G^{gp} = g\bar G ^p g^{-1}$.

We have the following result, whose proof is standard, from the fact that
$G_p \subset \bar G^p$.

\begin {lema} \label {aaab} $\mathcal U$ is an integrable distribution with integral
manifolds  $\bar G^q\cdot q$, $q\in M$.
\end {lema}

\begin {lema} \label {aaabb} $\mathcal U$ is a proper distribution of $M$ (or equivalently, since the bounded algebra $\mathfrak u ^p$ contains the isotropy algebra at $p$,  any bounded algebra is a proper subalgebra of $\mathfrak g$).
\end {lema}

\begin {proof}  Suppose on the contrary that $\mathcal{U}=TM$. Then for any $p\in M$, $\bar{ G}^p\cdot p = M$. Since
$\mathfrak u ^p = \mathrm {Lie}(\bar G^p)$, we would have that
$$\nabla _{\nu _p}X \subset \nu _p\, ,$$
for any Killing field induced by $\bar G^p$. Then, from Corollary
\ref {KostantdeRham}, $\nu $ is a parallel distribution. A contradiction.
\end {proof}

Summarizing, one has the following $G$-invariant distributions:
\begin {equation}\label {summarizing}
\{0\} \neq \nu \subsetneq  \nu ^{(1)} \subsetneq
 \nu ^{(2)}\subset  \mathcal U \subsetneq TM ,
\end {equation}
where $\nu$ and  $\nu ^{(1)}$ are  autoparallel and flat, and $\mathcal U$ is integrable and $G$-invariant.

\

 In the examples of last section we have that  $\nu ^{(2)} = \mathcal U$ and this distribution is not autoparallel, see Remark \ref{distriExamples}.

\begin {nota}\label {1321} Observe that (\ref {summarizing}) implies that the codimension $k$ of the nullity must be at least $3$. Moreover, if $k=3$ then
$\mathrm {codim}\,  \nu ^{(1)} = 2$,
$\nu ^{(2)}=  \mathcal U $, and  $\mathrm {codim}\, \mathcal U =1$.
\end {nota}

 \

 \begin {teo} \label {isodim} Let $M^n=G/H$ be a simply connected irreducible  homogeneous manifold with  nullity distribution $\nu$ of codimension
 $k$ ($G$ is not assumed to be connected). Then

 \begin {enumerate}

 \item [({\it i})]The representation $\rho$ of $H$ on $\nu _p^\perp$ is faithful ($p=[e]$).

 \item [({\it ii})] $\dim H \leq \frac 12 (k-2)(k-3)$. In particular, if $k= 3$, the isotropy is trivial (if $G$ is connected).

 \end {enumerate}
 \end {teo}

 \begin {proof} Assume that  $\ker (\rho)$ is a non-trivial normal subgroup of $H$.
Let $\mathbb V$ be the set of fixed vectors of $\ker (\rho)$.
Then $\nu _p ^\perp  \subset \mathbb V \subsetneq T_pM$. Then, by \cite {BCO}, Theorem 9.1.2,  $\mathbb V$ extends to an autoparallel $G$-invariant distribution $\mathcal D$ which contains the foliation $\nu  ^\perp$. Observe that for any $q \in M$,  $v\in \mathcal D_q,\,  w \in  \nu _q$ $R_{v,w} = 0$.
Then by Proposition 3.3 in \cite{DS}, $M$ splits. A contradiction. This proves (i).

  From (\ref {summarizing}), $H$ leaves any of the following subspaces invariant:
  $\nu _p \subsetneq \nu ^{(1)}_p\subsetneq \mathcal U _p\subsetneq T_pM$. Let us consider the representation $\rho$ of $H$ on the orthogonal complement $\nu_p^\perp$. Then there exist three non-trivial $H$-invariant mutually perpendicular subspaces $\mathbb V_1, \mathbb V _2, \mathbb V_3$ of dimensions $d_1$, $d_2$ and $d_3$ respectively such that

 $$\nu_p^\perp = \mathbb V_1 \oplus \mathbb V _2 \oplus \mathbb V_3, $$
 where $d_1\leq d_2\leq  d_3$, and $d_1+ d_2+  d_3 =k$.

 Then, by making use of part (i),

 \begin {eqnarray*}\dim H &=& \dim \rho (H) \leq \dim \mathrm {SO}(d_1) + \dim \mathrm {SO}(d_2) +
 \dim \mathrm {SO}(d_3) \\ &\leq& \dim \mathrm {SO}(k-2) = \frac 12 (k-2)(k-3),
\end{eqnarray*}
 where the last inequality is standard to show. This proves (ii).
\end {proof}

We finish this section by proving that $\nu ^{(1)}$ is a parallel distribution when restricted to any leaf of the bounded distribution.

%\color {red}
%\

\begin {lema}\label{parnu1} Let  $X,\,  U$ be  vector fields of $M$ such that $X$ lies in  $\nu ^{(1)}$ and $U$ lies in  $\mathcal U$. Then $\nabla _UX$ lies in $\nu ^{(1)}$ (and so  $\nu ^{(1)}$, restricted to any leaf of $\mathcal U$,  is a parallel and flat distribution).
\end{lema}

\begin {proof} Recall that $\nu ^{(1)} = \nu + \hat \nu$.
We may assume that either $X$ belongs to $\nu$ or $X$ belongs to $\hat \nu$.

({\it a}) Let $p\in M$ be arbitrary, let $X$ belong to $\nu$ and let $Z$ be a Killing field that belongs to  the bounded algebra $\mathfrak u ^p$, with $Z_p= U_p$. Since $\nu $ is $G$-invariant $[Z, X]$ lies in $\nu$.
Then $\nabla _{Z_p} X\in \nu _p$ if and only if
$\nabla _{X_p} Z \in \nu _p$, which follows from  the definition of the bounded algebra.

({\it b}) Let $X$ belong to $\hat \nu$. From the definition of $\hat \nu$, the fields of the form $\nabla _W Z$ span $\hat \nu$, where $Z\in \mathcal K ^G(M)$ and $W$ is a vector field of $M$ that lies in $\nu $. So we may assume that
$X= \nabla _W Z$.  Let $U$ be a vector field of $M$ that lies in $\mathcal U$.
From  (\ref {affK})
$\nabla ^2_{U,W} Z  =  R_{U,Z}W  = 0$. So,
$$0= \nabla ^2_{U,W} Z = \nabla _U\nabla _W Z -
\nabla _{\nabla _UW} Z.$$
Observe, from part (a),  that ${\nabla _UW}$ lies in $\nu$. Then
$\nabla _{\nabla _UW} Z$ lies in $\hat \nu$ and hence
$\nabla _U\nabla _W Z$ lies in $\hat \nu$.

 \end {proof}

%\color {black}

% Let us finish this section with a general fact.

% \begin {prop} \label {SSNCT}
% \end {prop}

 \section {Homogeneous spaces with  co-nullity $3$}

 Let $M^n=G/H$ be a simply connected Riemannian manifold with non-trivial nullity
  distribution $\nu$ of codimension $3$. Then, by Theorem \ref  {isodim},
  $H= \{e\}$ and so $M= G$, with a left invariant metric.  Then the autoparallel and
  $G$-invariant  distribution
  $ \nu ^{(1)}$ has codimension $2$, and the integrable $G$-invariant distribution
  $\mathcal U\supset \nu ^{(1)}$ has codimension $1$.

  Let $p\in M$ be fixed. Then there exist Lie subgroups $H_1\subset H_2$ of $G$ such that $H_1\cdot p \subset H_2\cdot p$ are the integral manifolds by $p$ of
  $ \nu ^{(1)}$ and $\mathcal U$, respectively (or, equivalently $\mathrm {Lie}
  (H_2) = \mathfrak u ^p$).

  \begin {lema} \label {setermina} We keep   the assumptions and notation of this section. Then

  \begin {enumerate}

  \item  $H_1\cdot p$ is isometric to $\mathbb R^{n-2}$.

 \item  $H_2\cdot p$ is  intrinsically flat.
 \
 \end {enumerate}

  \end {lema}

  \begin {proof}
That $H_1.p$ is flat was proved at the beginning of Section \ref{section5}.
By  Corollary \ref {STcor}, and part $(ii)$ of Proposition \ref{JacNul}, the Jacobi operator in any vector tangent to $\nu ^{(1)}$ is null.
So the proof of Lemma \ref {SC} also shows that the integral manifolds of $\nu ^{(1)}$ are simply connected.  This proves $(i)$.

Observe, since $k=3$, that  $\nu ^{(1)}$ is an autoparallel and flat sub-distribution of codimension $1$
of $\mathcal U$.
 From Lemma \ref {parnu1} it follows that $\nu ^{(1)}$, restricted to any leaf $S= H_2 \cdot p$ of $\mathcal U$ is a parallel and flat  distribution of codimension $1$. Then $S$ is flat, which proves (ii).
  \end {proof}

\vspace{.3cm}

\subsection{Excluding the Levi factors if $\mathbf {k=3}$ }
\label{sec61}
$\ $

Let $M$ be a homogeneous simply connected Riemannian manifold without Euclidean de Rham factor.
Assume that $M$ has a non-trivial nullity distribution of codimension 3. By Theorem \ref {isodim}, $M$ has no isotropy and so $M= G$, where $G=I(M)^o$ is endowed with a left invariant metric.

Let $\mathcal U$ be the bounded distribution,  which has codimension $1$ (see Remark \ref{1321}). The integral manifold by $p$ of $\mathcal U$ is given by
$H_2\cdot p$, where the Lie algebra of $H_2$ is the bounded algebra
 $\mathfrak u ^p$.
Since $G$ acts freely $\mathfrak u ^p$  has codimension $1$ in the Lie algebra
$\mathfrak g \simeq \mathcal K (M)$ of $G$.
Recall, from Lemma \ref {setermina} (ii), that $H_2\cdot p$ is intrinsically flat. Then, since there is no isotropy,
 $\mathfrak u^p$  is solvable (see Section \ref {homogeneous flat spaces}).

 Assume that the  Levi decomposition of $\mathfrak g$ has a non-trivial Levi factor. Namely,
 $$\mathfrak g =\mathfrak h \ltimes s, $$
 where $\mathfrak h$ is a semisimple Lie algebra and $s$ is the (solvable) radical of $\mathfrak g$.
 Observe that the intersection $\mathfrak h\cap \mathfrak u^p$ has codimension $1$ in $\mathfrak h$.
 Moreover, since $\mathfrak h$ is semisimple,  the projection of $\mathfrak u ^p$ to $\mathfrak h$ cannot be onto. This implies, since the codimension of $\mathfrak u ^p$ in $\mathfrak g$ is $1$, that
 \begin {equation}\label {Levi9}
 \mathfrak u ^p = \left( \mathfrak h\cap \mathfrak u^p \right) \ltimes \mathfrak s.
  \end {equation}
 Assume that $\mathfrak h = \mathfrak h _1 \oplus \cdots \oplus  \mathfrak h _r$
 is a direct sum of simple ideals. With the same arguments as before the projection $\mathfrak h'_i$ of
 $\mathfrak u ^p$ to $\mathfrak h _i$ has codimension $1$ in $\mathfrak h _i$ and coincides with $\mathfrak h_i\cap \mathfrak u^p $, $i=1, \cdots , r$.
We must have
 $\mathfrak u ^p =  (\mathfrak h'_1 \oplus \cdots \oplus \mathfrak h'_r) \ltimes \mathfrak s$, which implies,
  $r=1$, since $\mathfrak u ^p$ has codimension $1$.

If $\mathfrak h$ is of the compact type, then $\mathfrak h\cap \mathfrak u^p$ is solvable and so abelian.  Since $\mathfrak h$ is simple, this intersection must be properly contained in $\mathfrak h$ and so of codimension $1$.  Then $\mathfrak h\cap \mathfrak u^p$  is in the center of $\mathfrak h$ since each $\mathrm{ad}(x)$ is skew-symmetric w.r.t. the Killing form. A contradiction that shows that $\mathfrak h$  is of the non-compact type.

   Let $\mathfrak h = \mathfrak p \oplus \mathfrak k$ be the Cartan decomposition of
  $\mathfrak h$. Then $\mathfrak k$ has an $\mathrm {ad}_{\mathfrak k}$-invariant positive definite
  inner product $\langle \, , \,\rangle$. The intersection $ \mathfrak k\cap \mathfrak u ^p $ is solvable hence abelian and has codimension at most $1$. Since each $\mathrm{ad}(x)$, $x \in {\mathfrak k}$ is skew-symmetric w.r.t. $\langle \, , \,\rangle$ we conclude that $\mathfrak k$ is abelian. So, in any case, $\mathfrak h = \mathfrak {sl}_2$.

  We have shown, if $\mathfrak g$ is not solvable,  that
$$\mathfrak g = \mathfrak {sl}_2 \ltimes \mathfrak s, $$
where  $\mathfrak s$ is the  radical of $\mathfrak g$.
Let
$\mathfrak b := \mathfrak {sl}_2\cap \mathfrak u ^p$, which is solvable since $\mathfrak u ^p$ is so.
As previously observed,
$$\mathfrak u ^p= \mathfrak b \ltimes \mathfrak s. $$

Since $\mathfrak b $ is solvable,  there exist,  as it is well-known, a basis $A, B, C$ of
$\mathfrak {sl} _2$, such that $A, B$ span $\mathfrak b$ and
\begin {equation} \label {sl2tri}
[A, B] = 2B, \ \ \ [A, C] = -2C,  \ \ \ [B, C] = A
\end {equation}
i.e. $A, B, C$ is a so-called $\mathfrak {sl}_2$-triple.

We will identify any element $v\in \mathfrak g$ with the Killing field $q\mapsto v.q$ of $M$. This identification is a Lie algebra anti-isomorphism. With this identification, after replacing $A$ by $-A$, we have the same relations of (\ref {sl2tri}) for $A, B, C$.

\

\begin {lema}\label {nablaAB}
 We keep the assumptions and notation of this section. Then
\begin {enumerate}
\item   $ \nabla _{\nu _p} B = \{0\}$.
\item $\nabla _{\nu _p} A = \{0\}$.
\item $\nabla _{\nu _p} C = \{0\}$.
\end {enumerate}
\end {lema}

\begin {proof}

From equation \ref{curvatureformula} we have that

$$ 2\nabla B = \nabla [A, B] = R_{A, B} - [\nabla A , \nabla B] . $$
Then
\begin {eqnarray}\label {12453}
2(\nabla B)_{\vert \nu _p} = (\nabla  [A, B])_{\vert \nu _p}&=& - [\nabla A , \nabla B]_{\vert \nu _p}\\
&=& -[(\nabla A)_{\vert \nu _p} , (\nabla B)_{\vert \nu _p}].
\end {eqnarray}
where last equality follows from the fact that $A, B \in \mathfrak u ^p$ and so $\nabla _{\nu _p}A , \nabla _{\nu _p}B \subset \nu _p$.
But the skew-symmetric endomorphism $ -[(\nabla A)_{\vert \nu _p} , (\nabla B)_{\vert \nu _p}]$ of
$\nu _p$ is perpendicular to  $(\nabla B)_{\vert \nu _p}$ (with the usual inner product).
Then $\nabla _{\nu _p}B = \{0\}$ which proves (i).
\vspace {.2cm}

From \ref{curvatureformula}  we have that
\begin {eqnarray}\label {aa994}
  \nabla _{\nu _p}A &=& \nabla_ {\nu _p} [B, C]
= - (\nabla _{\nabla  _{\nu _p} C} B  - \nabla _{\nabla _{\nu _p}B} C) \\
&=& - \nabla _{\nabla  _{\nu _p}C} B,
\end {eqnarray}
where the last equality is due to (i).

On the one hand,  $ - \nabla _{\nabla C _{\nu _p} }B$ must be perpendicular to
$\ker (\nabla B)_p \supset \nu _p$. On the other hand, $\nabla _{\nu _p}A \subset \nu _p $. Then, from \ref {aa994}, we obtain that $\nabla _{\nu _p}A   = \{0\}$, and so (ii).
\vspace {.2cm}

Recall  that  $[A, C] = -2C$, and let $v\in \nu _p$ be arbitrary. Then,  by equation \ref {curvatureformula}
\begin {eqnarray}
-2\nabla _{v}C &=& \nabla _{v}[A, C]= - (\nabla _{\nabla  _{v} C} A  - \nabla _{\nabla _{v}A} C) \\
& = &
 - \nabla _{\nabla  _{v} C} A,
 \end {eqnarray}
 where the last equality is due to (ii).
 Since $(\nabla A)_p$ is skew-symmetric the last term of the above equality is perpendicular to
 $\nabla  _{v} C$. But the first term of this equality is proportional to $\nabla  _{v} C$.
 Then $\nabla  _{v} C= 0$, which proves (iii).
\end {proof}

\
Lemma \ref {nablaAB} implies that $C$ belongs to the bounded algebra $\mathfrak u ^p$.
But $\mathfrak g$ is linearly spanned by $C$ and $\mathfrak u ^p$. Then
$$\mathfrak g  =\mathfrak u ^p.$$
 This contradicts Lemma \ref {aaabb}.
Then $\mathfrak g$ has no Levi factor and so we obtain the following result:

\begin {teo}\label {Gsolvable} Let $M=G/H$ be a simply connected homogeneous Riemannian manifold without Euclidean de Rham factor. Assume that the nullity distribution of $M$ is non-trivial and of codimension $k=3$. Then $H= \{e\}$ and $G$ is solvable.
\end {teo}

\medskip
\section {The leaves of $\nu$ are closed and $\nu ^\perp$ is completely non-integrable} \label {LTA}
\smallskip

\begin {lema}\label {EDF} Let $M= G/H$ be a (non-simply connected) homogeneous Riemannian manifold and let $\mathcal {D}^0$ be the parallel distribution of $M$ associated to its  local Euclidean de Rham factor of $M$  ($G$ connected).  Let $\bar F (p)$ be the closure of a (maximal)  integral manifold $F(p)$ by $p$ of $\mathcal {D}^0$. Then there is a closed abelian   normal  subgroup $A$  of $ I(M)$ (not depending on $p$ and non-necessarily contained in $G$) such that
$\bar F (p)  = A\cdot p$, for all $p\in M$. In particular, $\bar F (p)$ is a
flat (embedded) homogeneous  submanifold of $M$.
\end {lema}

\begin {proof}

Let $\tilde M$ be the universal cover of $M$  and write it as
$\tilde M = \mathbb R^k\times M_1$, where $\mathbb R^k$ is the Euclidean
de Rham factor. Let $\tilde G \subset I(\tilde {M})$ be the (connected) lift of $G$. Let
$\Gamma $ be the deck transformations of $\tilde {M}$. Then $\Gamma$ commutes with
$\tilde G$. Let $\Gamma ^0$ be the image of the projection of $\Gamma$ to $I(\mathbb R^k)$. Let $\tilde G^0$ and $\tilde G^1$ be the images of the projections of $\tilde G$ into
$I(\mathbb R^k)$ and $I(M_1)$, respectively. Then $\Gamma ^0$ commutes with the transitive group $\tilde G^0$ of isometries. Then the elements of $\Gamma ^0$ are translations (see Section \ref {homogeneous flat spaces}).
Let $\tilde {T} \simeq \mathbb R ^k$ be group of translations of $\mathbb R ^k$.
Then $\tilde {T} \times \tilde G^1$ acts transitively on $\tilde M$ and commutes with $\Gamma$. Then $T \times \tilde G^1$  projects to a transitive group of isometries of $M$. Let $T  \subset I(M)$ be the projection of $\tilde T$. Observe,  since $\tilde {T}$ is a normal subgroup of $I(\tilde M)$, that  $T$ is a normal subgroup of $I(M)$.
Let $A$ be the closure of $T$ in $I(M)$. Then $A$ is a normal subgroup of $I(M)$ and
$\bar {F}(p) = A\cdot p $.
\end {proof}

\begin {nota}\label {REDF}  Observe that in  the assumptions and notation of Lemma \ref {EDF} it follows that the family of closures of the integral manifolds of
$\mathcal D^0$ is an $I(M)$-invariant foliation of $M$.
\end {nota}

\medskip

 Let $M= G/H$ be a simply connected irreducible homogeneous Riemannian manifold with a non-trivial nullity distribution $\nu$, where $G=I(M)^o$ (the identity component of the isometry group of $M$).
 Let, for $p\in M$,
 $E^p$ be the Lie subgroup of $G$ that leaves invariant  the integral manifold
 $N(p)$ of $\nu$ by $p$ (see Section \ref{homogeo}). Observe that since $\nu$ is $G$-invariant,  $G_p= (E^p)_p$. We may assume that $p=[e]$ so that  $G_p=H$.
 Then, since $N(p)$ and $H$ are connected,  $E^p$ is connected.
 Let  $\bar E^p$ be the closure of $E^p$ in $G=I(M)^o$ and let
$\bar N (p)$  be the closure of  $N(p)$ in $M$.
Then
$$ \bar N (p) = \bar E^p\cdot p .$$

Since $\nu _p\subset T_p\bar N (p)$, then $\bar \nu = \nu _{\vert \bar N (p)}$ is a distribution
of $\bar N(p)$. Since
$E^p$ is a normal subgroup of $\bar E^p$, then the integral manifolds of the autoparallel distribution $\bar \nu$ are given by
\begin {equation}\label {766}
 E^p\cdot x  \ \ \ \ \ (x\in \bar N(p)).
\end {equation}

\begin {lema}\label {nubarN} $\bar \nu $ is contained in the nullity distribution of
$\bar N(p)$.
\end {lema}

\begin {proof} Let us first show that $\bar \nu$ is in the nullity of the second fundamental form $\alpha$ of $\bar N(p)$.
Let $X$ be a Killing field of $M$ induced by $\bar E^p$ and let $\gamma _v (t) =
\mathrm {Exp} (tu)p$ be a homogeneous geodesic in $N(p)$, $u\in \mathrm {Lie}
(E^p)$ with $u.p=v$.
The proof of Proposition \ref {existstransv} shows that there is transvection on the direction of $\nabla _vX$. Moreover, from the construction of such a transvection,
one has that $\nabla _vX\in T_p\bar N(p)$). Since the Killing field $X$, induced by
$\bar E^p$, is arbitrary, we conclude that $\nu_p$ belongs to the nullity of
$\alpha$ at $p$ (and the same is true for any $q\in \bar N(p)$).
Then, from the Gauss equation, one obtains that $\bar \nu$ is contained in the nullity distribution of $\bar N(p)$.
\end {proof}

\begin {nota}  Observe that the proof of Lemma \ref {nubarN} applies also to  prove  the following fact: let  $G'$ be a Lie subgroup of $G$ that contains $E^p$ and let $S = G'\cdot p$. Then $\nu _{\vert S}$ is contained in the nullity distribution of the Riemannian manifold $S$.
\end {nota}

\begin {teo}\label {MMLPQTP} Let $M= G/H$ be a simply connected irreducible homogeneous  Riemannian manifold with non-trivial nullity distribution $\nu$, where $G =I(M)^o$.Then  any (maximal) integral manifold of $\nu$ is a closed (embedded) submanifold of $M$.
\end {teo}

\begin {proof}
Applying Proposition \ref {KillTanNu} (and the comment below it) any Killing field, induced by
$\bar E^p$, that lies in the nullity  of
 $\bar N (p)$ must be tangent to $\mathcal D ^0$, the distribution associated to the
 local Euclidean de Rham factor of $\bar N(p)$.  Since the Killing fields induced by $E^p$ lie in $\nu$, which is included in the nullity of $\bar N (p)$, then
 $N(p)$ is included in the integral manifold $F(p)$ of $\mathcal D ^0$.
 By Lemma \ref {EDF} the closure $\bar F (p)$ is a flat embedded  submanifold of
 $M$.  Observe that $\bar F (p) = T^r\times \mathbb R ^s$, since it is flat and homogeneous (see Section \ref {homogeneous flat spaces}).

 By Remark \ref {REDF} there is a maximal (connected) Lie subgroup $S$ of $\bar E^p$ that contains $E^p$ and  such that
 $S\cdot p = \bar F (p)$. Then $E^p$ is a normal subgroup of $S$.
 Then the orbits $N(y) = E^p\cdot y$  are parallel totally geodesic  submanifolds of
 $\bar F (p)$, $y\in \bar F (p)$
 (see Lemma \ref {wwqq}).
 Let us consider the intersection $T_pN(p)\cap T_pT^r$. If this intersection is different from $\{0\}$, then $N(p)$ has a bounded geodesic $\gamma _v$. Then any Killing field of $M$ is bounded along $\gamma _v$ and hence, by Lemma \ref
 {ABC}, parallel along $\gamma _v$. So, for any $X\in \mathcal K^G(M)$, $\mathbb Rv$ is invariant under $(\nabla X)_p$. Then, by Corollary \ref {KostantdeRham}, $M$ splits off the direction of $v$. A contradiction.
So $T_pN(p)\cap T_pT^r = \{0\}$. In this case, as it is not hard to see,
$N(p)$ is a closed submanifold of $\bar F (p)$. Then $N(p)$ is closed or,  equivalently, $E^p$ is a closed subgroup of $I(M)$.
\end {proof}

\begin {lema}\label {wwqq}
Let $G$ be a connected subgroup of $I(\mathbb R ^n)$ which acts transitively on
$\mathbb R ^n$. Let $G'$ be connected normal subgroup of $G$.
Then the orbits of $G'$ are parallel affine subspaces of $\mathbb R ^n$.
\end {lema}

\begin {proof} It is well-known that any (connected) Lie subgroup of $I(\mathbb R ^n)$ has a totally geodesic orbit (see e.g. Theorem 3.5, pg. 100 in \cite{AVS}). Then, since $G'$ is a normal subgroup of $G$, all orbits of $G'$ are affine subspaces of
$\mathbb R ^n$.  Let $d$ be the distance between
the affine subspaces $G'\cdot x$ and $G'\cdot y$. We may assume that $d(x,y) = d$.
So for any $g'x\in G'\cdot x$,   $d(g'x, g'y) = d$. Then any point of $G'\cdot y$ is at a distance $d$ to some point in $G'\cdot x$. Then the affine subspace $G'\cdot x$ must be parallel to the affine subspace $G'\cdot y$.
\end {proof}

\smallskip

\subsection {The affine bundle and the connection given by $\mathbf {\nu ^\perp}$}\
\smallskip

We keep the assumptions and notation of this section. From Theorem  \ref{MMLPQTP}
we have that $N(p)$ is closed, or equivalently, $E^p$ is a closed subgroup
of $G=I(M)^o$ ($p=[e]$). Then $M=G/H$ is the total space of a fiber bundle over
$B= G/E^p$ with standard fiber $E^p/H= N(p)\simeq \mathbb R^k$ (with a Euclidean affine structure, see Lemma \ref {SC}). The projection of $M$ onto $B$ will be denoted by $\pi$.
Observe that  $B$ is the quotient space $M/\mathcal N$ of $M$
 by the leaves of the nullity foliation $\mathcal N = \{N(q): q\in M\}$.
 Since the elements of $\mathcal N$ are, in a natural way, Euclidean affine spaces, one has that $M$ is an Euclidean affine bundle (and so an affine combination of local  sections  of $M$ is a local section). Observe that $G$ leaves $\mathcal N$ invariant and, for any
$g\in G$, $g$ is an isometry between  $\pi ^{-1}(\pi (q)) = N(q)$ and $\pi ^{-1}(\pi (gq))=N(gq) $. Moreover,  Lemma \ref {KillTanNu} implies that $G$ acts almost effectively on  $M/\mathcal{N}=B$.

Let us consider the natural affine connection on $M\overset {\pi}{\to}B$ given by the distribution $\nu ^\perp$. In fact, a  perpendicular variation of totally geodesic manifolds is by isometries. So, the local horizontal lift of curves in $B$ gives rise to  local isometries between the involved fibers. From this particular situation,  it is well-known, and standard to show, that any piece-wise differentiable curve
$c:[0,1]\to B$ can be lifted to a (unique) horizontal curve $\tilde c_u:[0,1]\to M$
with $\tilde c_u (0) = u$, for any
$u\in \ \pi ^{-1}(c(0))$.  Then there is a well defined parallel transport
$\tau _c: \pi ^{-1}(c(0)) \to \pi ^{-1}(c(1))$, which is an isometry, given by
$\tau _c (u) = \tilde c _u (1)$. Then $\nu ^\perp$ is an affine connection.

For each $b\in B$, let $\Phi (b)\subset I(\pi ^{-1}(b))$ denote the holonomy group of
$\nu ^\perp$ at $p$ (holonomy groups are conjugated by parallel transport). Note that $B$ is simply connected, since $M$ is simply connected and the fibers are connected. Then the holonomy  groups $\Phi (b)$  are connected.

Let us consider, for $q\in M$, the holonomy subbundle
$\mathrm {Hol}(q)$.
Namely, $\mathrm {Hol}(q)$ consists of all the elements of $M$ that can be reached
from $q$ by a horizontal curve. The holonomy subbundles foliate $M$
. Moreover, any holonomy subbundle intersects any given fiber $\pi ^{-1}(b)$ in an orbit of the holonomy group $\Phi (b)$.

 The holonomy subbundles, despite what happens in a principal bundle,
 may have different  dimensions depending on the dimensions of the orbits of the holonomy group. But in our case the holonomy subbundles have all the same dimension (and so their tangent spaces define a smooth distribution), since, for any $g\in G$, $u\in M$,
\begin{equation} \label {holhol}
\mathrm {Hol}(gq) = g\mathrm {Hol}(q) \text { \ \ \ and \ \ \ }
 \Phi (\pi (g q)) =g \Phi (\pi (q))g^{-1}.
  \end{equation}

 In particular, if $g\in E^q$ then  $\Phi (\pi (q)) =g \Phi (\pi (q))g^{-1}$.

 This implies that $\tilde E^q := \{g_{\vert N(q)} : g\in E^q\}$ is included in the normalizer of $\Phi (\pi (q))$ in $I(N(q))$.
 Then $L = \Phi (\pi (q)) . \tilde{E}^q$ is a Lie group of isometries, which is transitive
 on $N(q)\simeq \mathbb R ^k$, and $\Phi (\pi (q))$ is a normal subgroup of $L$.
 Then, by Lemma \ref {wwqq}, we have that:

\centerline {  \ ($*$) \    {\it The orbits of $\Phi (\pi (q))$ are parallel affine subspaces of $N(q)$.}}
\

Observe that the above property implies that $\Phi (\pi (q))$ {\it acts polarly}  on
$N(q)$.

 \medskip

 Let $\mathcal {Y}$ be the distribution of $M$ defined by the normal spaces of
 the holonomy subbundles. Namely,
 \begin {equation}\label {abcrrr}  \mathcal {Y}_q = (T_q\mathrm {Hol}(q))^\perp \subset \nu _q.
 \end {equation}

 Then $\mathcal {Y}\subset \nu$ and, from ($*$),
 $\mathcal {Y}_{\vert N(q)}$ is a parallel (i.e. constant) foliation of $N(q)$, for all $q\in M$ (which is perpendicular to the holonomy orbits). Moreover, since all the orbits of $\Phi (\pi (q))$ are principal orbits:
 \smallskip

 \centerline {  ($**$) \ \   {\it  any $w\in \mathcal {Y}_q$ is a fixed vector of the isotropy $\Phi (\pi (q))_q$.}}

 Let $p\in M$ be fixed and let $\xi \in \mathcal {Y}_p$. Then $\xi$ induces a normal vector field of $\mathrm {Hol}(p)$ in the following way: if $q\in \mathrm {Hol}(p)$, choose  $\tilde c_p : [0,1] \to M$ be a horizontal piece-wise differentiable curve
 with $\tilde c_p(0)=p$, $\tilde c_p(1)=q$. Let $c=\pi \circ \tilde c _p$, then
  define
  $$\tilde \xi (q)= \mathrm {d}\tau _c (\xi).$$

  From ($**$) one obtains that $\tilde \xi$ is well defined (and it is standard to show that it is smooth). Let us show that $\tilde \xi$ is a parallel normal vector field.
  Observe that

 $$T_p\mathrm {Hol}(p) = \nu ^\perp_p
 \oplus T_p(\Phi (\pi (p))\cdot p).$$
From the construction of $\tilde \xi$, taking into account that $\Phi (\pi (p))$ acts polarly on $N(p)$, one obtains that
 $$\nabla ^\perp _v\tilde \xi= 0,  \text {\ \ \ \
 for all } v\in {T_p(\Phi (\pi (p))\cdot p)}.$$

 Let now $u\in \nu ^\perp _p$ and let $\tilde c _p (s)$ be a horizontal curve with
 $\tilde c^\prime _p (0)= u$  and let $c(t) =\pi (\tilde c_p(t))$. Let $\gamma _\xi (t)$ be a geodesic (i.e. a line) in $N(p)$ with $\gamma ^\prime_\xi (0) =\xi$.
 Let $\tau ^s$ be the parallel transport along $c(s)$, form $0$ to $s$.
 Let us consider $$f(s,t) = \tau ^s (\gamma _\xi (t)), $$
 which is variation of geodesics that are tangent to  $\nu$.
 Then
 $$J_ \xi (t):=\tfrac {\partial f }{\partial s}_ {\vert (0, t)}$$
 is a Jacobi field along the geodesic
 $\gamma _{\xi}(t)$.

 Note  that $s\mapsto f(s,t_0)$ is a horizontal curve,
 i.e.  tangent to $\nu ^\perp$, and so $J^\xi (t)$ is a horizontal field along
 $\gamma _\xi (t)$.

Since $N(p)$ is totally geodesic in $M$,   $\nu ^\perp _{\vert N(p)}$ is a parallel subbundle of the pull-back of $TM$, via the inclusion of
$N(p)$ in $M$.

Then
\begin {equation} \label {4544}
J_\xi ^\prime(t) := \tfrac {\mathrm {D}\, }{ \mathrm {d}t}J_\xi (t),
\end {equation}
as well as $J_{\xi}(t)$, are horizontal fields along $\gamma _\xi (t)$.

Observe that
\begin {eqnarray*}\label {2233}
J_\xi ^\prime(0) &=& \tfrac {\mathrm {D}\, }{ \partial t}_{\vert (0,0)}
\tfrac {\partial \,}{\partial s}f (s,t)=
\tfrac {\mathrm {D}\, }{ \partial s}_{\vert (0,0)}
\tfrac {\partial \,}{\partial t}f (s,t)  \\
&=& \tfrac {\mathrm {D}\, }{ \mathrm {d} s}_{\vert 0}\mathrm {d}\tau ^s (\xi)
= \tfrac {\mathrm {D}\, }{ \mathrm {d} s}_{\vert 0}\tilde \xi (\tilde c_p(s)) \\
&=& \nabla _u\tilde \xi \in \nu ^\perp _p.
\end{eqnarray*}

Since $u\in \nu ^\perp_p$ is arbitrary, then   $\tilde \xi$ is also $\nabla ^\perp$-parallel in the horizontal directions and so
$$\nabla ^\perp \tilde \xi =0. $$

Observe, since $\xi\in (T_p(\mathrm {Hol}(p)))^\perp$ is arbitrary, that:

\smallskip
\centerline {\ ($***$) \ \ {\it The normal bundle of
$\mathrm {Hol}(p)$ is globally flat. }}
\smallskip

Let $A$ be the shape operator of $\mathrm {Hol}(p)$. It is standard to show, since $\nabla ^\perp\tilde \xi =0$, that $J_\xi ^\prime(0) = -A_\xi (u)$. So, the Jacobi field $J_\xi (t)$ has the following initial conditions:
\begin {eqnarray*}\label {7545}
\ J_\xi (0)&=& u, \\
J_\xi ^\prime (0)&=&  -A_\xi (u).
\end {eqnarray*}
 Since $J_\xi ^\prime (0)\in \nu^\perp_p$, this shows that
any shape operator of $\mathrm {Hol}(p)$ leaves $\nu ^\perp _p$ invariant.

Note that $J_\xi (t) = \tilde c _{\gamma _\xi (t)}^\prime (0)$, where
$\tilde c _{\gamma _\xi} (t)$ is the horizontal lift of $c(t)=\pi (\tilde c _p(t))$.
Then
\begin {equation}\label {3976} \mathrm {d}\pi (J_\xi (t) ) = c^\prime (0)= \mathrm {d}\pi (u).
\end {equation}

Assume that $\lambda \neq 0$ is an eigenvalue of  $A_\xi$ and let $u\neq 0$ be an eigenvector of $A_\xi$ associated to $\lambda$. Since $\gamma _\xi (t)$ lies in
$N(p)$, then
$$J_\xi (t) = (1-t\lambda)\hat u (t),$$
where $\hat u (t)$ is the parallel transport of $u$ along $\gamma _\xi (t)$.
Then $J(1/\lambda) =0$, which contradicts (\ref {3976}). This shows that
$$A_{\xi \vert \nu ^\perp_p}=0= A_{\tilde \xi (p) \vert \nu ^\perp_p}.$$

From its  construction, $\tilde \xi$ is constant along the holonomy orbit
$\Phi (\pi (p))\cdot p$ (see ($*$) and ($**$)).
Then
$$A_{\xi \vert T_p(\Phi (\pi (p))\cdot p)}=0$$
and so
$$A_\xi =0.$$
Since $p$ and $\xi \in (T_p(\mathrm {Hol}(p)))^\perp =\mathcal Y _p$ are arbitrary, we conclude that any holonomy subbundle  $\mathrm {Hol}(q)$ is a totally geodesic submanifold of $M$.
Observe that the distribution $\mathcal Y$ is autoparallel, since $\mathcal Y\subset \nu$ and is parallel inside the leaves of the nullity
(see the paragraph below (\ref
{abcrrr})).

Since the holonomy subbundles are totally geodesic and its perpendicular distribution $\mathcal Y$ is autoparallel, we conclude that $\mathcal Y$ is a parallel distribution. This is a contradiction, since $M$ is irreducible, unless $\mathcal Y =0$. Then $\mathrm {Hol}(p) =M$ and so $\Phi (p)$ is transitive on $N(p)$.

 By summarizing our main results in the section we obtain:

\begin {teo}\label{transihol} Let $M=G/H$ be a simply connected  irreducible homogeneous Riemannian manifold with a non-trivial nullity distribution $\nu$.
Then the quotient space $B$ of $M$ by the leaves of the nullity is a manifold and
$M$ is a Euclidean affine fiber bundle over $B$, with standard fiber isometric to $\mathbb R^k$. Moreover, $\nu ^\perp$ defines a metric  affine connection on $M$  with a transitive holonomy group (and so $\nu ^\perp$ is completely non-integrable).
\end {teo}

For Euclidean or spherical submanifolds, the  complete non-integrability of the distribution perpendicular to the relative nullity (i.e. the nullity of the second fundamental form) was proved in
\cite {Vi}  (this is not true in hyperbolic space).

\section{The proof of the main results}\label{demostraciones}
The main theorems stated in the introduction were proved throughout the paper or are direct consequences of previous results. We sum them up here.
\medskip

\noindent
{\it Proof of Theorem A.} The first part of (1) was proved throughout Section \ref{section5} and concluding with equation (\ref {summarizing}).
 The fact that the integral manifolds of $\nu$ are simply connected was proved in Lemma \ref {SC}.
By  Corollary \ref {STcor}, and part $(ii)$ of Proposition \ref{JacNul}, the Jacobi operator in any vector tangent to $\nu ^{(1)}$ is null.
So the proof of Lemma \ref {SC} also shows that the integral manifolds of $\nu ^{(1)}$ are simply connected.

Part (2). The existence of an adapted transvection $Y$, see Definition \ref{adapted}, with $Y_p=v \notin \nu_p$ was proved in Proposition \ref{existstransv}. The fact that the Jacobi operator $R_{,v}v$ is null was proved in Theorem \ref {ST} and stated in Corollary \ref{STcor}. Theorem \ref {ST} also states that  $[Y,[Y, \mathcal K(M)]] = 0$ and that $Y$ does not belong to the center of $\mathcal K ^G(M)$.
For an arbitrary $v\in \hat {\nu} _p$ the existence of a transvection $Y$ with the stated properties  follows from Proposition  \ref {JacNul}, since the adapted transvections at $p$ span $\hat \nu _p$.

Part (3) is Theorem \ref {isodim}.

The first part of (4) is Theorem \ref {Gsolvable}.
In Section \ref {EX33} we construct non-trivial examples, in any dimension, with $k=3$ and $G$ not unimodular. For the non-unimodularity  see Remark \ref{Milnor}.  \qed
\vspace {.2cm}

\noindent
{\it Proof of Theorem B.} It follows from Theorem \ref {MMLPQTP} and Theorem
\ref {transihol}. \qed

\vspace {.2cm}

\noindent
{\it Proof of Proposition C.} Part (a) is Corollary \ref{corsemisimple}. If $M$ is compact see Proposition \ref{compact2} for a direct proof. Part (b) is Corollary \ref {cornilpotent}. \qed

%-------------------
%-------------------

\color {black}

\section{Examples with nullity of codimension $3$ and co-index of symmetry $2$}
\label{EX33}

In this section we construct examples of irreducible Riemannian homogeneous spaces with nullity of codimension 3 and co-index of symmetry 2 in any dimension greater or equal to $4$. As explained in the introduction such examples are optimal since neither the nullity can be greater nor the co-index of symmetry can be smaller due to Reggiani's Theorem \cite{R}.

Let $G = \mathbb{R}^{d} \rtimes \mathbb{R}  $ be the semidirect product of the abelian groups where $\mathbb{R}$ acts on $\mathbb{R}^d$ as
$\mathrm {Exp} (sA) = \mathrm {e}^{sA}$, $s \in \mathbb{R}$ and $A$ is defined as
\[ A = \begin{bmatrix}a M_{d-1} & a e_1 \\
          -a e_1^t & a
  \end{bmatrix} \]
where $M_{d-1} = (m_{ij})$ is $(d-1) \times (d-1)$ skew-symmetric with $m_{ij} = 1$ for $i < j$, $e_1$ is the first canonical column of $\mathbb{R}^{d-1}$, $e_1^t$ its transpose,  and the constant $a$ is chosen so $1 = trace(A A^t)$ ,  i.e. $a^2 = \frac{1}{3 + (n-2)(n-3)}$,  where $n=d+1$. \\

 Observe that $G$ is simply connected. We can regard $G$ as  the connected Lie subgroup of $\mathrm{GL}(d+1,\mathbb{R})$ whose Lie algebra is generated by the following $d+1$ matrices:
\[ E_i := \begin{bmatrix} 0 & e_i \\
          0 & 0
  \end{bmatrix} , \hspace{1cm} \A := \begin{bmatrix} A & 0\\
          0 & 0
  \end{bmatrix} \]
where $i=1,\dots,d$ and $e_1,\cdots,e_d$ are the canonical columns of $\mathbb{R}^d$.

Let $g$ be the left invariant metric on $G \subset \mathrm{GL}(d+1,\mathbb{R})$   whose restriction to the tangent space at the identity $T_{\ed}G$ is given by
$$g(X,Y) = \mathrm{trace}(X Y^t), $$ 
where $Y^t$ indicates the transpose matrix.
Observe that $E_1,\cdots,E_d,\A$ is an orthonormal base of $T_{\ed}G$.

 Let $\tilde E_1,\cdots,\tilde E_d,\tilde {\A}$ be the right invariant vector fields of $G$ such that
$$(\tilde E_1)_{\ed }= E_1,\cdots,(\tilde E_d)_{\ed } = E_d,(\tilde {\A})_{\ed } = \A.$$
Observe that the flow associated to any of these fields  is given by the left multiplication by $\mathrm {Exp}(tu)$, where $u$ is the evaluation of the field at $e$. So $\tilde E_1,\cdots,\tilde E_d,\tilde {\A}$ are Killing fields of $G$, since their associated flows are by isometries.  For the sake of simplifying the notation we will identify, when there is no possible confusion,
\begin{equation}\label {identif}\tilde E_1 \simeq  E_1,\cdots, \tilde E_d \simeq  E_d,  \tilde {\A}\simeq  \A.\end{equation}

\begin{lema}\label{nablas} Let $\nabla$ be the Levi-Civita connection of $(G,g)$ and $R$ its curvature tensor.
Then  the following identities  hold at $\ed \in G$:

\begin{itemize}
\item[i)] $\nabla E_1 = \nabla E_2 = \cdots = \nabla E_{d-1} = 0$, i.e. $E_1,\cdots,E_{d-1}$ are transvections at $\ed \in G$.

\item[ii)] $\nabla E_d =
\begin{bmatrix}
0 & 0 & 0 \\
0 & 0 & a \\
0 &-a & 0 \\
\end{bmatrix} $.

\item[iii)] $\nabla \A =
\begin{bmatrix}
-a M_{d-1} & -a e_1 & 0  \\
a e_1^t & 0 & 0\\
0 & 0 & 0
\end{bmatrix}$
\end{itemize}
where the matrices  have to be regarded  as (skew-symmetric) linear endomorphisms of $T_{\ed}G$  in the orthonormal basis  $E_1, \cdots , E_d, \, \A$.
\end{lema}

\begin{proof} The proof is a computation by using equation (\ref{fundamentalequation})  and the fact that bracket at $\ed$ of any two right invariant fields is minus the bracket at $\ed$ of the corresponding left invariant fields (see (\ref{laspe2})).
To show ({\it i}) notice that equation (\ref{fundamentalequation}) gives $\langle \nabla_{E_i} E_j, E_k \rangle = 0$ for any $i,j,k \in {1,\cdots,d}$. Equation (\ref{fundamentalequation}) gives
\[ 2 \langle \nabla_{\A} E_j, E_k \rangle = \langle [\A , E_j],E_k \rangle + \langle [\A , E_k],E_j\rangle.\]
 So, for $ 1 \leq j \leq d-1$, we get
\[ 2 \langle \nabla_{\A} E_j, E_k \rangle = \langle Ae_j ,e_k \rangle + \langle Ae_k ,e_j \rangle = a m_{jk} + a m_{kj} = 0 \]
which shows ({\it i}).\\
To show ({\it ii}) observe that from the definition of $A$, if $1<k<d$:
\[ 2 \langle \nabla_{\A} E_d, E_k \rangle = \langle Ae_d ,e_k \rangle + \langle Ae_k ,e_d \rangle = \langle a e_1 + a e_d ,e_k \rangle + \langle Ae_k ,e_d \rangle = 0 =  2 \langle C^aE_d , E_k\rangle \]
and
 \[ 2 \langle \nabla_{\A} E_d, E_d \rangle = \langle Ae_d ,e_d \rangle + \langle Ae_d ,e_d \rangle = 2a =
 2\langle C^a\A, E_d\rangle\, ,
  \]
 where $C^a = \begin{bmatrix}
0 & 0 & 0 \\
0 & 0 & a \\
0 &-a & 0 \\
\end{bmatrix} $.
 This shows ({\it ii}).

Finally, \[ 2 \langle \nabla_{E_i} {\A} , E_j \rangle = \langle [E_i,\A],E_j \rangle + \langle [\A,E_j],E_i \rangle = 2 \langle \nabla_{E_i} {\A} , E_j \rangle = -\langle A e_i,e_j \rangle + \langle A e_j,e_i \rangle.  \]
So,  if $i,j \in \{1,\cdots,d-1\}$, we have 
\[ 2\langle \nabla_{E_i} {\A} , E_j \rangle =  -a m_{ij}+a m_{ji} = -2a m_{ij} = 2\langle D^aE_i, E_j\rangle,  \]
where $ D^a=\begin{bmatrix}
-a M_{d-1} & -a e_1 & 0  \\
a e_1^t & 0 & 0\\
0 & 0 & 0
\end{bmatrix}$

Now for $1<j<d$ we get
\begin {eqnarray*} 2\langle \nabla_{E_d} {\A} , E_j \rangle &=& -\langle A e_d,e_j \rangle + \langle A e_j,e_d \rangle = -\langle a e_1 + a e_d,e_j \rangle + \langle A e_j,e_d \rangle = 0 \\
&=&  2\langle D^aE_d, E_j\rangle.
\end {eqnarray*}

If $j=1$, then
\[2\langle \nabla_{E_d} {\A} , E_1 \rangle = -\langle A e_d,e_1 \rangle + \langle A e_1,e_d \rangle = -2a     =  2\langle D^aE_d, E_1\rangle.\]
The last case to check is
\[2\langle \nabla_{\A} {\A} , E_j \rangle = 0   = 2\langle D^aA, E_j\rangle.\]
 Thus, since $\nabla \A$ is skew-symmetric,  we have proved ({\it iii}).
\end{proof}

\begin{lema}\label {ABAB}  The nullity $\nu_{\ed}$ of $R$ at $\ed \in G$ is generated by $E_2,\cdots,E_{d-1}$.
\end{lema}
\begin{proof} We are going to compute $R_{X, E_i}$ for $i=1,\cdots,d-1$ by using formula (\ref{curvatureformula}).  Let $X$ be an arbitrary Killing field of $(G,g)$.  We regard  $E_1, \cdots , E_d, \A$, from the identification \ref {identif}, as Killing fields of $(G,g)$.
We have that $$R_{X, E_i} = \nabla [X,E_i],$$
since $\nabla E_i = 0$ by   Lemma \ref{nablas},   part ({\it i}). So
\[R_{X, E_i} = \nabla [\sum_{j=1}^{d} \langle X,E_j \rangle E_j + \langle X,\A \rangle \A ,E_i] = \langle X,\A \rangle \nabla [\A ,E_i] = \langle X,\A \rangle R_{\A , E_i} \, .\]
By formula (\ref{curvatureformula}) $R_{\A , E_i}=\nabla [\A , E_i]$ then
\[ R_{\A , E_i} = \nabla \left( \sum_{k=1}^{d-1} a m_{ki} E_k  -a \delta_{1i} E_d \right) \]
and by Lemma \ref{nablas} we get $R_{\A ,E_i} = -a \delta_{1i} \nabla E_d$. Then $E_2,\cdots,E_{d-1}$ belongs to the nullity $\nu_{\ed}$ of $R$.
Since $R_{\A ,E_1} = -a \nabla E_d$, Lemma  \ref{nablas} ({\it ii}) implies that \[ \nu_{\ed} = span\{E_2,\cdots,E_d\} \]
because  of  $R_{\A E_1} \neq 0$ and $\mathrm {ker}(\nabla E_d) = span\{E_1,\cdots,E_{d-1}\}$.
\end{proof}

\begin{nota}\label{distriExamples} Observe that in these examples the Killing vector field $Y$ defined as
\[ Y:= E_1 - \sum_{k=3}^{d-1} E_k \] is a transvection as  that given by  part ({\it ii}) of Theorem {A}. Indeed, by part ({\it iii})  of Lemma \ref{nablas}, $Y =  \frac 1a\nabla_{E_2} \A$.  Hence $Y_e \in {\hat \nu}_e$ and $Y_e \notin \nu_e$, and the Jacobi operator $R_{\cdot, Y}Y$ is null:
\[ R_{X, Y}Y = R_{X , E_1}E_1 = \langle X,\A \rangle R_{\A , E_1}E_1 = \langle X,\A \rangle (-a) \nabla_{E_1} E_d = 0 \, .\]
Actually, in our examples, $\nu^{(1)}_{\ed}=span\{E_1, E_2,\cdots,E_{d-1}\}$ and
$$\nu^{(2)}_{\ed}={\mathcal{U}}_{\ed} = span\{E_1, E_2,\cdots,E_d\}  .$$
\end{nota}

\begin{lema}\label{nonablainv} There is no (non-trivial)  $\nabla \A$-invariant proper subspace of $\nu_{\ed}$.
\end{lema}
\begin{proof}  From  part ({\it iii}) of Lemma \ref{nablas} and
Lemma \ref {ABAB},
any non-trivial  invariant subspace of $\nabla \A$ in $\nu_{\ed}$  determines a non-trivial invariant subspace of the matrix $M_{d-1}$  which is contained in the subspace of $\mathbb{R}^{d-1}$ generated by the vectors $e_2,\cdots,e_{d-1}$. But this is not possible by the lemma in  Appendix.
\end{proof}

To show that $(G,g)$ is an irreducible Riemannian manifold we use the following  result:

\begin{lema}\label {prepre} Let $M=G/H$ be a homogeneous Riemannian manifold whose nullity distribution has codimension 3.
If $M$ is locally reducible,  then the flat local de Rham factor is non-trivial.
\end{lema}
\begin{proof} If there is  no local flat factor, then there must be an irreducible local factor whose nullity has codimension 1.   This implies that all sectional curvatures are zero (since any plane of the tangent space does intersect the nullity subspace). Hence this factor must be flat.   A contradiction.
\end{proof}

 \begin{teo}\label{examples} For each $n \geq 4$ the $n$-dimensional simply connected homogeneous Riemannian manifold 
  $(G,g)$, constructed in this section, is irreducible. Its nullity distribution has codimension 3 and its co-index of symmetry is $2$. Moreover its Ricci tensor has four eigenvalues: zero with multiplicity $n-3$, $-a^2$ and $a^2 \left( \frac{-1 \pm \sqrt{5}}{2}\right)$. So $(G,g)$ has sectional curvatures of both signs and its scalar curvature is $-2a^2 = \frac{-2}{3 + (n-2)(n-3)}$.
\end{teo}
\begin{proof}
That the index of symmetry $i_\mathfrak{s}(G)$ is $n-2$ follows from $i)$ in Lemma \ref{nablas}.
 Let us show that  our examples are irreducible Riemannian manifolds. In fact, if $(G,g)$ is reducible,  by   Lemma \ref {prepre} there is a non-trivial flat de Rham factor ${E}$ of $(G,g)$.
Observe that the parallel $G$-invariant distribution $\mathcal E$ of $(G,g)$, associated to the factor $E$,  is contained in the nullity distribution $\nu$.
Moreover $[\A , \mathcal E] \subset \mathcal E$ since $\mathcal E$ is is invariant by isometries.
Let $Z$ be a parallel vector field of $(G,g)$ that lies in $\mathcal E$.
Then $([\A , Z])_{\ed} = (\nabla_Z \A)_{\ed} - (\nabla_{\A} Z)_{\ed} = (\nabla_Z \A)_{\ed} \in \mathcal E_{\ed} \subset \nu_{\ed}$.
Since the parallel vector fields \color {black} of $(G,g)$ generate $\mathcal E$, it follows that
$\nabla \A$ leaves invariant the non-trivial  subspace $ \mathcal E_{\ed}$ of $\nu_{\ed}$. This contradicts Lemma \ref{nonablainv}.
A direct computation using Lemma \ref{nablas} and formula (\ref{curvatureformula}) shows that the  Ricci tensor restricted to $\nu_{\ed}^{\perp} = \mathrm{span}\{E_1,E_d,\A \}$ is given by the matrix $\begin{bmatrix} 0 & a^2 & 0 \\ a^2 & -a^2 & 0 \\ 0 & 0 & -a^2 \end{bmatrix}$  (with respect to the basis $E_1,E_d,\A$).
\end{proof}

\begin{nota} For $d=3$, by changing the matrix $A$, it is possible to show the existence
of 1-parameter family $(G_{\lambda},g_{\lambda})$ of 4-dimensional
 non-homothetic irreducible homogeneous
metrics with nullity of dimension $1$.
\end{nota}

\begin{nota} With the same ideas,  by modifying the matrix $A$, it is possible to construct examples with $k>3$.
\end{nota}

\begin{nota}\label{Milnor} Observe that $\mathrm{trace}(\mathrm{ad}(\A)) \neq 0$ i.e. our solvable groups are not unimodular hence they do not admit finite volume quotients \cite[Remark, Lemma 6.2.]{Mi76}.
\end{nota}

\section{Appendix: invariant subspaces of the skew-symmetric matrix $M$}\label{skewM}

Let $M = (m_{ij})$ be the real $d \times d$, $(d>2)$, skew-symmetric matrix with $m_{ij} = 1$ if $i<j$.
Let $\mathbb{W} \subset \mathbb{R}^d$ be the subspace generated by the canonical vectors $e_2,\cdots,e_d$ i.e. the orthogonal complement of the first canonical vector $e_1$.
The goal of this appendix is to prove the following:

\begin{lema}  There are no $M$-invariant, non-trivial,  subspaces of $\mathbb R^d$ which are contained in $\mathbb{W}$.
\end{lema}

\begin{proof}
By contradiction assume that there is a non-trivial $M$-invariant subspace $\mathbb{U} \subset \mathbb{W}$.
We have to consider two cases $dim(\mathbb{U})=1$ or $dim(\mathbb{U})=2$.

Case $dim(\mathbb{U})=1$. Since $M$ is skew-symmetric we have that $\mathbb{U} \subset \ker(M)$. Let $(a_1,\cdots,a_d) \neq 0 \in \mathbb{U} \subset \ker(M)$.
Then
\begin{equation} \label{defM} M.a = (\sum_{k=2}^{d}a_k , \cdots , -\sum_{k=1}^{i-1}a_k + \sum_{k=i+1}^d a_k , \cdots,  -\sum_{k=1}^{d-1}a_k) \, .
\end{equation}
So subtracting two consecutive components we obtain:
\begin{center}
$\begin{cases} 0 = a_2 + a_1 , \\ 0 = a_3 + a_2 \\ \, \, \, \, \, \, \vdots \\ 0 = a_d + a_{d-1} .
\end{cases}$
\end{center}
Since $a \in \mathbb{U} \subset \mathbb{W}$ we have $a_1 = 0$ hence $a=0$. A contradiction.

Case $dim(\mathbb{U})=2$. In this case $\mathbb{U}$ is spanned by two vectors $a, b:=M.a$. Since $M$ is skew-symmetric there is $r \neq 0 \in \mathbb{R}$ such that $M.b = r.a$. We can assume $r \neq 1$ since otherwise the vector $a+b$ is $M$-invariant which was excluded in Case 1.
By using equation \ref{defM} we get \begin{center}
$\begin{cases} b_1 - b_2 = a_1 + a_2 \\ ra_1 - ra_2 = b_1 + b_2
\end{cases}$
\end{center}
By using $a,b \in \mathbb{U} \subset \mathbb{W}$ we get that $a_1=b_1=0$ and so $\begin{cases}  - b_2 = a_2 \\  - ra_2 =  b_2
\end{cases}$ hence $a_2=b_2=0$ due to $r\neq1$.
Now equation \ref{defM} gives us for $i=1,\cdots,d-1$
\begin{center}
$\begin{cases} b_i - b_{i+1} = a_i + a_{i+1} \\ ra_i - ra_{i+1} = b_i + b_{i+1}
\end{cases}$ \,.
\end{center}
So if we assume, as inductive hypothesis, that $a_1=a_2=\cdots=a_i = b_1=b_2=\cdots=b_i = 0 $ we obtain for $a_{i+1},b_{i+1}$:

\begin{center}
$\begin{cases}  - b_{i+1} =  a_{i+1} \\ - ra_{i+1} =  b_{i+1}
\end{cases}$ \,.
\end{center}
Hence $a_{i+1},b_{i+1} = 0$ due to $r \neq 1$. Then $a=b=0$.  A contradiction.

\end{proof}

%\subsection*{Acknowledgments}
%The authors thank
%Leandro Cagliero for useful discussions.

\vspace{1cm}

\noindent
\begin{tabular}{l|cl|cl}
Antonio J. Di Scala & & Carlos E. Olmos  & & Francisco Vittone\\
\footnotesize Dipartimento di Scienze Matematiche&  & \footnotesize FaMAF, CIEM-Conicet & &\footnotesize DM, ECEN, FCEIA - Conicet\\
\footnotesize Politecnico di Torino &  &\footnotesize Universidad Nac. de C\'ordoba& &\footnotesize Universidad Nac. de Rosario \\
 \footnotesize Corso Duca degli Abruzzi, 24& &\footnotesize Ciudad Universitaria  & &\footnotesize Av. Pellegrini 250\\
 \footnotesize 10129, Torino, Italy  & & \footnotesize 5000, C\'ordoba, Argentina & &\footnotesize 2000, Rosario, Argentina \\
\footnotesize{antonio.discala@polito.it} & & \footnotesize{olmos@famaf.unr.edu.ar} & & \footnotesize{vittone@fceia.unr.edu.ar}\\
\footnotesize{http://calvino.polito.it/$\sim$adiscala}& &
 & &\footnotesize{www.fceia.unr.edu.ar/$\sim$vittone}
\end{tabular}

\end {document}